\documentclass[11pt]{article}
\usepackage{srcltx}
\usepackage{epsfig}
\usepackage{url}
\usepackage[colorlinks=true]{hyperref}
\def\bZ{{\mathbf{Z}}}

\def\VI{{\hbox{\rm VI}}}

\def\epsilonsad{\epsilon_{\hbox{\tiny\rm sad}}}
\def\epsilonVI{\epsilon_{\hbox{\tiny\rm VI}}}
\def\epsilonNash{\epsilon_{\hbox{\tiny\rm Nash}}}
\def\Res{{\hbox{\rm Res}}}

\def\SV{{\hbox{\rm SadVal}}}
\def\cL{{\cal L}}
\usepackage{psfrag}
\usepackage{graphicx}
\usepackage[latin1]{inputenc}
\usepackage{subfigure}
\usepackage{hyperref}
\usepackage{mathtools}
\usepackage{amsfonts,amsmath,amssymb,amsthm}
\usepackage{color}
\usepackage{booktabs}

\usepackage{float,enumitem}
\setlist{itemsep=3pt,parsep=0pt,topsep=3pt}

\floatstyle{ruled}
\newfloat{algorithm}{tbp}{loa}
\floatname{algorithm}{Algorithm}

\DeclareMathOperator{\rank}{rank}

\newtheorem{thm}{Theorem}

\newtheorem{lemma}[thm]{Lemma}

\newcommand{\cP}{{\cal P}}
\def\qed{\ \hfill$\square$\par\smallskip}

\def\cV{{\cal V}}

\def\Opt{{\mathop{\hbox{Opt}}}}

\def\rank{{\mathop{\hbox{\rm Rank}}}}

\def\Min{\mathop{\rm Min}}
\def\Max{\mathop{\rm Max}}

\def\cY{{\cal Y}}

\newcommand{\cW}{{\cal W}}

\newcommand{\be}{\begin{eqnarray}}
\newcommand{\ee}[1]{\label{#1}\end{eqnarray}}

\newcommand{\ese}{\end{eqnarray*}}
\newcommand{\bse}{\begin{eqnarray*}}

\newtheorem{proposition}{Proposition}
\newtheorem{remark}{Remark}

\def\Argmin{\mathop{\hbox{\rm Argmin$\,$}}}

\def\Argmax{\mathop{\hbox{\rm Argmax$\,$}}}

\def\cU{{\cal U}}
\def\cA{{\cal A}}
\def\cS{{\cal S}}
\def\cD{{\cal D}}
\def\cH{{\cal H}}
\def\cI{{\cal I}}

\def\Card{{\hbox{\rm Card}}}

\def\bR{{\mathbf{R}}}

\def\cB{{\cal B}}

\def\cX{{\cal X}}
\def\cJ{{\cal J}}
\definecolor{darkmagenta}{RGB}{125,38,205}

\begin{document}
\title{Decomposition Techniques for Bilinear Saddle Point Problems and Variational Inequalities with Affine Monotone Operators on Domains Given by Linear Minimization Oracles}
\author{Bruce Cox\thanks{US Air Force} \and Anatoli Juditsky\thanks{LJK,
Universit\'e J. Fourier, B.P. 53, 38041 Grenoble
Cedex 9, France, {\tt Anatoli.Juditsky@imag.fr}}
\and Arkadi Nemirovski\thanks{Georgia Institute
 of Technology, Atlanta, Georgia
30332, USA, {\tt nemirovs@isye.gatech.edu}\newline
Research of
the third author was supported by the NSF grants CMMI-1232623, CCF-1415498, CMMI-1262063.}
}
\maketitle
\begin{abstract}
 The majority of First Order methods for large-scale convex-concave saddle point problems and variational inequalities with monotone operators are {\sl proximal} algorithms which at every iteration need to minimize over problem's domain $X$ the sum of a linear form and a strongly convex function. To make such an algorithm practical, $X$ should be {\sl proximal-friendly} -- admit a strongly convex function with easy to minimize linear perturbations. As a byproduct, $X$ admits a computationally cheap Linear Minimization Oracle (LMO) capable to minimize over $X$ linear forms. There are, however, important situations where a cheap LMO indeed is available, but $X$ is not proximal-friendly, which motivates search for algorithms based solely on LMO. For smooth convex minimization, there exists a classical LMO-based algorithm -- Conditional Gradient.  In contrast, known to us LMO-based techniques for other problems with convex structure (nonsmooth convex minimization, convex-concave saddle point problems, even as simple as bilinear ones, and variational inequalities with monotone operators, even as simple as affine)  are quite recent and utilize common approach based on Fenchel-type representations of the associated objectives/vector fields. The goal of this paper is to develop an alternative (and seemingly much simpler) {\sl decomposition} LMO-based techniques for {\sl bilinear} saddle point problems and for variational inequalities with {\sl affine} monotone operators. \end{abstract}
 \section{Introduction}
 This paper is a follow-up to our paper \cite{FTVI} and, same as its predecessor, is motivated by the desire to develop first order algorithms for solving  convex-concave saddle point problem (or variational inequality with  monotone operator) on a convex domain $X$ represented by {\sl Linear Minimization Oracle} (LMO) capable to minimize over $X$, at a reasonably low cost, any linear function. ``LMO-representability'' of a convex domain $X$ is an essentially weaker assumption than ``proximal friendliness'' of $X$ (possibility to minimize over $X$, at a reasonably low cost, any linear perturbation of a properly selected strongly convex function) underlying the vast majority of known first order algorithms.  There are important applications giving rise to LMO-represented domains which are {\sl not} proximal friendly, most notably
 \begin{itemize}
 \item nuclear norm balls arising in low rank matrix recovery and in Semidefinite optimization; here LMO reduces to approximating the leading pair of singular vectors of a matrix, while all known proximal algorithms  require much costly computationally full singular value decomposition,
 \item total variation balls arising in image reconstruction; here LMO reduces to solving a specific flow problem \cite{withZaid}, while a proximal algorithm needs to solve a much more computationally demanding linearly constrained convex quadratic program,
 \item some combinatorial polytopes.
\end{itemize}
The needs of these applications inspire the current burst of activity in developing LMO-based optimization techniques. In its major part, this activity was focused on smooth (or Lasso-type smooth regularized) Convex Minimization over LMO-represented domains, where the classical Conditional Gradient algorithm of Frank \& Wolfe \cite{FW} and its modifications are applicable (see, e.g.,
\cite{DemRub,Dunn,Freund,GarHazan,HDPDM.11,withZaid,Jag1.16,Jag2.17,PsheD} and references therein). LMO-based techniques for large-scale Nonsmooth Convex Minimization (NCM), convex-concave Saddle Point problems (SP), even bilinear ones, and  Variational Inequalities (VI) with monotone operators, even affine ones, where no classical optimization methods work, have been developed only recently. To the best of our knowledge, the related results reduce to
LMO-based techniques for large-scale NCM based on Nesterov's smoothing \cite{Argyriou2014,Pierucci2014,Tewari2011,Ying2012,DualSubgradient,Lan1}. An alternative approach to NCM, based on Fenchel-type representations of convex functions and processing the induced by these representations problems dual to the problem of interest, was  developed in \cite{DualSubgradient} and was further extended in  \cite{FTVI} to convex-concave SP's and VI's with monotone operators.
The goal of this paper is to develop an alternative to \cite{FTVI} {\sl  decomposition-based} approach to solving convex-concave SP's and
monotone VI's on LMO-represented domains. In the nutshell, this approach is extremely simple, and it makes sense to present an {\sl informal} outline of it, in the SP case,  right here.
 \begin{quote}
 Given convex compact sets $X_1,X_2,Y_1,Y_2$ in Euclidean spaces, consider a convex-concave saddle point ``master'' problem
$$
\min\limits_{[x_1;x_2]\in X_1\times X_2}\max_{[y_1;y_2]\in Y_1\times Y_2}\Phi(x_1,x_2;y_1,y_2)
$$
along with two ``induced'' problems
$$
\begin{array}{ll}
(P)&\min_{x_1\in X_1}\max_{y_1\in Y_1}\left[\phi(x_1,y_1):=\min_{x_2\in X_2}\max_{y_2\in Y_2}\Phi(x_1,x_2;y_1,y_2)\right]
\\
(D)&\min_{x_2\in X_2}\max_{y_2\in Y_2}\left[\psi(x_2,y_2):=\min_{x_1\in X_1}\max_{y_1\in Y_1}\Phi(x_1,x_2;y_1,y_2)\right]
\\
\end{array}
$$
It is easily seen that $(P)$ and $(D)$ are convex-concave problems and a good approximate solution to the master problem
induces straightforwardly equally good approximate solutions to $(P)$ and to $(D)$. More importantly, it turns out that when solving one of the induced problems, say, $(P)$, by an ``intelligent,'' in certain precise sense, algorithm, information acquired in course of building an $\epsilon$-solution yields straightforwardly an $\epsilon$-solution to the master problem, and thus yields an $\epsilon$-solution to the other induced problem, in our case, to $(D)$.\\
$\quad$
Now imagine that  we  want to solve a convex-concave SP problem which ``as is'' is too complicated for the standard solution techniques (e.g., problem's domain is not proximal-friendly, or is of huge dimension). Our proposed course of actions is to make the problem of interest the problem $(D)$ stemming from a master problem built in a way which ensures that the associated problem $(P)$ is amenable to an  ``intelligent'' solution algorithm $\cB$. After such a master problem is built, we solve $(P)$ within a desired accuracy $\epsilon$ by $\cB$ and use the acquired information to build an $\epsilon$-solution to the problem of interest.
\end{quote}
As we shall see, our decomposition approach can, in principle, handle general convex-concave SP's and affine VI's. Our emphasis in this paper is, however, on {\sl bilinear} SP's and on VI's with {\sl affine} operators -- the cases which, on one hand, are of primary importance in numerous applications, and, on the other hand, are the cases where our  approach  is easy to implement and where this approach seems to be more flexible and much simpler than the machinery of Fenchel-type representations developed in \cite{FTVI} (and in fact even covers this machinery, see section \ref{sect:dominance}).
\par
The rest of this paper is organized as follows. In sections \ref{sect:SPdecomposition} and \ref{sect:monotoneVI} we present our decomposition-based approach to SP problems, resp., VI's with monotone operators, with emphasis on utilizing the approach to handle {\sl bilinear} SP's, resp., {\sl affine} VI's, on {\sl LMO-represented} domains. We illustrate our constructions by applying them to  Colonel Blotto type matrix game (section \ref{Blotto}) and Nash Equilibrium with pairwise interactions (section \ref{sect:skewsymm}); in both these illustrations, decomposition allows to overcome difficulties coming from potentially huge ambient dimensions of the problems.\\
Proofs missing in the main body of the paper are relegated to Appendix.

\section{Decomposition of Convex-Concave Saddle Point Problems}\label{sect:SPdecomposition}
\subsection{Situation}\label{s:sit}
In this section, we focus on the situation as follows. Given are
\begin{enumerate}
\item convex compact sets $X_i$ in Euclidean spaces $\cX_i$ and convex compact sets $Y_i$ in Euclidean spaces $\cY_i$, $i=1,2$;
\item convex compact sets $X$, $Y$ such that
$$
X\subset X_1\times X_2\subset\cX:=\cX_1\times \cX_2,\,Y\subset Y_1\times Y_2\subset\cY:=\cY_1\times \cY_2,
$$
      such that the projections of $X$ onto $\cX_i$ are the sets $X_i$, and projections of $Y$ onto $\cY_i$ are the sets $Y_i$, $i=1,2$. For $x_1\in X_1$, we set $X_2[x_1]=\{x_2:[x_1;x_2]\in X\}\subset X_2$, and for $y_1\in Y_1$ we set $Y_2[y_1]=\{y_2:[y_1;y_2]\in Y\}\subset Y_2$. Similarly,
$$
X_1[x_2]=\{x_1:[x_1;x_2]\in X\},\,x_2\in X_2,\ \hbox{and}\ Y_1[y_2]=\{y_2:[y_1;y_2]\in Y\},\, y_2\in Y_2;
$$
\item Lipschitz continuous function
\begin{equation}\label{Phi}
\Phi(x=[x_1;x_2];y=[y_1;y_2]):X\times Y\to\bR
\end{equation}
which is convex in $x\in X$, and concave in $y\in Y$.
\end{enumerate}
We call the outlined situation a {\sl direct product} one, when $X=X_1\times X_2$ and $Y=Y_1\times Y_2$.
\subsection{Induced convex-concave functions}\label{s:iccf}
We associate with $\Phi$ {\sl primal and dual induced} functions:
$$
\begin{array}{rcl}
\phi(x_1,y_1)&:=&\min\limits_{x_2\in X_2[x_1]}\max\limits_{y_2\in Y_2[y_1]}\Phi(x_1,x_2;y_1,y_2)=\max\limits_{y_2\in Y_2[y_1]}\min\limits_{x_2\in X_2[x_1]}\Phi(x_1,x_2;y_1,y_2):X_1\times Y_1\to\bR,\\
\psi(x_2,y_2)&:=&\min\limits_{x_1\in X_1[x_2]}\max\limits_{y_1\in Y_1[y_2]}\Phi(x_1,x_2;y_1,y_2)=\max\limits_{y_1\in Y_1[y_2]}\min\limits_{x_1\in X_1[x_2]}\Phi(x_1,x_2;y_1,y_2):X_2\times Y_2\to\bR.\\
\end{array}
$$
(the equalities are due to the convexity-concavity and continuity of $\Phi$ and convexity and compactness of $X_i[\cdot]$ and $Y_i[\cdot]$).\par
Recall that a Lipschitz continuous convex-concave function $\theta(u,v):U\times V\to\bR$ with convex compact $U,V$ gives rise to the primal and dual problems
$$
\begin{array}{rcll}
\Opt(P[\theta,U,V])&=&\min\limits_{u\in U}\left[\overline{\theta}(u):=\max\limits_{v\in V}\theta(u,v)\right]&\\
\Opt(D[\theta,U,V])&=&\max\limits_{v\in V}\left[\underline{\theta}(v):=\min\limits_{u\in U}\theta(u,v)\right]&\\
\end{array}
$$
with equal optimal values:
$$
\SV(\theta,U,V):=\Opt(P[\theta,U,V)]=\Opt(D[\theta,U,V]),
$$
same as gives rise to {\sl saddle point residual}
$$
\epsilonsad([u;v]|\theta,U,V)=\overline{\theta}(u)-\underline{\theta}(v)=\left[\overline{\theta}(u)-\Opt(P[\theta,U,V])\right]+
\left[\Opt(D[\theta,U,V])-\underline{\theta}(v)\right].
$$
\begin{lemma}\label{lem1} $\phi$ and $\psi$ are convex-concave on their domains, are lower (upper)  semicontinuous in their ``convex''  (``concave'') arguments, and are Lipschitz continuous in the direct product case. Besides this,  it holds
\begin{equation}\label{eq9}
\SV(\phi,X_1,Y_1)=\SV(\Phi,X,Y)=\SV(\psi,X_2,Y_2),
\end{equation}
and whenever $\bar{x}=[\bar{x}_1;\bar{x}_2]\in X$ and $\bar{y}=[\bar{y}_1;\bar{y}_2]\in Y$, one has
\begin{equation}\label{eq10}
\epsilonsad([\bar{x}_1;\bar{y}_1]|\phi,X_1,Y_1)\leq\epsilonsad([\bar{x};\bar{y}]|\Phi,X,Y),\,\,
\epsilonsad([\bar{x}_2;\bar{y}_2]|\psi,X_2,Y_2)\leq\epsilonsad([\bar{x};\bar{y}]|\Phi,X,Y).\\
\end{equation}
\end{lemma}
\paragraph{The strategy} for solving SP problems we intend to develop is as follows:
\begin{enumerate}
\item We represent the SP problem of interest as the {\sl dual SP problem}
$$\min_{x_2\in X_2}\max_{y_2\in Y_2}\psi(x_2,y_2)\eqno{(D)}$$ induced by {\sl master SP problem}
    $$
    \min_{[x_1;x_2]\in X}\max_{[y_1;y_2]\in Y}\Phi(x_1,x_2;y_1,y_2)\eqno{(M)}
    $$
The master SP problem is built in such a way that the associated {\sl primal SP problem}
$$
\min_{x_1\in X_1}\max_{y_1\in Y_1}\phi(x_1,y_1)\eqno{(P)}
$$
admits First Order oracle and can be solved by a traditional First Order method (e.g., a proximal one).
\item We solve $(P)$ to a desired accuracy by First Order algorithm producing {\sl accuracy certificates} \cite{NOR} and use these certificates to recover approximate solution of required accuracy to the problem of interest.
\end{enumerate}
We shall see that the outlined strategy (originating from \cite{CoxPhD}\footnote{in hindsight, a special case of this strategy was used  in \cite{Gol96,Gol97}.}) can be easily implemented when the problem of interest is a bilinear SP on the direct product of two LMO-represented domains.
 \subsection{Regular sub- and supergradients}\label{s:reg}
 Implementing the outlined strategy requires some ``agreement'' between the first order information of the master and the induced SP's, and this is the issue we address now.
 \par
Given $\bar{x}_1\in X_1,\bar{y}_1\in Y_1$, let $\bar{x}_2\in X_2[\bar{x}_1]$ and $\bar{y}_2\in Y_2[\bar{y}_1]$ form a saddle point of the function
$\Phi(\bar{x}_1,x_2;\bar{y}_1,y_2)$ ($\min$ in $x_2\in X_2[\bar{x}_1]$, $\max$ in $y_2\in Y_2[\bar{y}_1]$). In this situation we say that $(\bar{x}=[\bar{x}_1;\bar{x}_2],\bar{y}=[\bar{y}_1;\bar{y}_2])$ belongs to the {\sl saddle point frontier} of $\Phi$, and we denote this frontier by $\cS$.
 Let now $\bar{z}=(\bar{x}=[\bar{x}_1;\bar{x}_2],\bar{y}=[\bar{y}_1;\bar{y}_2])\in\cS$, so that the function $\Phi(\bar{x}_1,x_2;\bar{y}_1,\bar{y}_2)$ attains its minimum over $x_2\in X_2[\bar{x}_1]$ at $\bar{x}_2$, and the function  $\Phi(\bar{x}_1,\bar{x}_2;\bar{y}_1,y_2)$  attains its maximum over $y_2\in Y_2[\bar{y}_1]$ at $\bar{y}_2$. Consider a subgradient
 $G$ of $\Phi(\cdot;\bar{y}_1,\bar{y}_2)$ taken at $\bar{x}$ along $X$: $G\in \partial_x\Phi(\bar{x};\bar{y})$. We say that $G$ is a {\sl regular} subgradient of $\Phi$ at $\bar{z}$,  if for some $g\in E_1$ it holds
$$
\forall x=[x_1;x_2]\in X: \langle G,x-\bar{x}\rangle \geq \langle g,x_1-\bar{x}_1\rangle;
$$
every $g$ satisfying this relation is called {\sl compatible} with $G$.
Similarly,  we say that  a supergradient $H$ of $\Phi(\bar{x};\cdot)$, taken at $\bar{y}$ along $Y$ is a regular supergradient of $\Phi$ at $\bar{z}$,  if for some $h\in F_1$ it holds
$$
\forall y=[y_1;y_2]\in Y: \langle H,y-\bar{y}\rangle \leq \langle h,y_1-\bar{y}_1\rangle,
$$
and every $h$ satisfying this relation will be called {\sl compatible} with $H$.
 \begin{remark}\label{rem1} Let $X=X_1\times X_2$, $Y=Y_1\times Y_2$, meaning that we are in the direct product case.
 If $\Phi(x;\bar{y})$ is differentiable in $x$ at $x=\bar{x}$, the partial gradient
 $\nabla_x\Phi(\bar{x};\bar{y})$ is a regular subgradient of $\Phi$ at $(\bar{x},\bar{y})$, and $\nabla_{x_1}\Phi(\bar{x};\bar{y})$ is compatible with this subgradient:
$$
\begin{array}{l}
\forall x=[x_1;x_2]\in X_1\times X_2:\\
 \langle \nabla_x\Phi(\bar{x};\bar{y}),x-\bar{x}\rangle
=\langle\nabla_{x_1}\Phi(\bar{x};\bar{y}),x_1-\bar{x}_1\rangle+
  \underbrace{\langle \nabla_{x_2}\Phi(\bar{x};\bar{y}),x_2-\bar{x}_2\rangle}_{\geq0} \geq
  \langle \nabla_{x_1}\Phi(\bar{x};\bar{y}),x_1-\bar{x}_1\rangle.
  \end{array}
  $$
 Similarly, if $\Phi(\bar{x};y)$ is differentiable in $y$ at $y=\bar{y}$, then the partial gradient
 $\nabla_y\Phi(\bar{x};\bar{y})$ is a regular supergradient of $\Phi$ at $(\bar{x},\bar{y})$, and $\nabla_{y_1}\Phi(\bar{x};\bar{y})$ is compatible with this supergradient.
 \end{remark}

 \begin{lemma}\label{lem2} In the situation of section \ref{s:sit}, let $\bar{z}=(\bar{x}=[\bar{x}_1;\bar{x}_2],\bar{y}=[\bar{y}_1;\bar{y}_2])\in\cS$, let $G$ be a regular subgradient
 of $\Phi$ at $\bar{z}$ and let $g$ be compatible with $G$. Let also $H$ be a regular supergradient of $\Phi$ at $\bar{z}$, and $h$ be compatible with $H$. Then $g$ is a subgradient in $x_1$, taken at $(\bar{x}_1,\bar{y}_1)$ along $X_1$, of the induced function $\phi$, and $h$ is a supergradient in $y_1$, taken at $(\bar{x}_1,\bar{y}_1)$ along $Y_1$, of the induced function $\phi$:
 $$
 \begin{array}{ll}
 (a)&\phi(x_1,\bar{y}_1)\geq\phi(\bar{x}_1;\bar{y}_1)+\langle g,x_1-\bar{x}_1\rangle,\\
 (b)&\phi(\bar{x}_1,y_1)\leq\phi(\bar{x}_1;\bar{y}_1)+\langle h,y_1-\bar{y}_1\rangle.\\
 \end{array}
 $$
for all $x_1\in X_1$, $y_1\in Y_1$.
\end{lemma}

\paragraph{Regular sub- and supergradient fields of induced functions.} In the sequel, we  say that $\phi^\prime_{x_1}(x_1,y_1)$,
$\phi^\prime_{y_1}(x_1,y_1)$ are {\sl regular} sub- and supergradient fields of $\phi$, if for every $(x_1,y_1)\in X_1\times Y_1$ and properly selected $\bar{x}_2$, $\bar{y}_2$ such that the point
$\bar{z}=(\bar{x}=[x_1;\bar{x}_2],\bar{y}=[y_1;\bar{y}_2])$ is on the SP frontier of $\Phi$, $\phi^\prime_{x_1}(x_1,y_1)$, $\phi^\prime_{y_1}(x_1,y_1)$ are the sub- and supergradients of $\phi$ induced, via Lemma \ref{lem2}, by regular sub- and supergradients of $\Phi$ at $\bar{z}$. Invoking Remark \ref{rem1}, we arrive at the following observation:
\begin{remark}\label{rem2} Let $X=X_1\times X_2$, $Y=Y_1\times Y_2$, meaning that we are in the direct product case.
If $\Phi$ is differentiable in $x$ and in $y$,  then regular sub- and supergradient fields of $\phi$ can be built as follows: given $(x_1,y_1)\in X_1\times Y_1$, we find  $\bar{x}_2$, $\bar{y}_2$ such that the point $\bar{z}=(\bar{x}=[x_1;\bar{x}_2],\bar{y}=[y_1;\bar{y}_2])$ is on the SP frontier of $\Phi$, and set
\begin{equation}\label{eq1}
\phi^\prime_{x_1}(x_1,y_1)=\nabla_{x_1}\Phi(x_1,\bar{x}_2;y_1,\bar{y}_2),\,\,\phi^\prime_{y_1}(x_1,y_1)=\nabla_{y_1}\Phi(x_1,\bar{x}_2;y_1,\bar{y}_2).
\end{equation}
\end{remark}
\subsubsection{Existence of regular sub- and supergradients}
The notion of regular subgradient deals with $\Phi$ as a function of $[x_1;x_2]\in X$ only, the $y$-argument being fixed, so that the existence/description questions related to regular subgradient deal in fact with a Lipschitz continuous convex function on $X$. And of course the questions about existence/description of regular supergradients reduce straightforwardly to existence/decription of regular subgradients (by swapping the roles of $x$'s and $y$'s and passing from $\Phi$ to $-\Phi$). Thus, as far as existence and description of regular sub- and supergradients is concerned, it suffices to consider the situation where
\begin{itemize}
\item $\Psi(x_1,x_2)$ is a Lipschitz continuous convex function on $X$,
\item $\bar{x}_1\in X_1$, and $\bar{x}_2\in X_2[\bar{x}_1]$ is a minimizer of $\Psi(\bar{x}_1,x_2)$ over $x_2\in X_2[\bar{x}_1]$.
\end{itemize}
What we need to understand, is when a subgradient $G$ of $\Psi$  taken at $\bar{x}=[\bar{x}_1;\bar{x}_2]$ along $X$ and some $g$ satisfy the relation
\begin{equation}\label{regularity}
\langle G,[x_1;x_2]-\bar{x}\rangle \geq \langle g,x_1-\bar{x}_1\rangle\,\forall x=[x_1;x_2]\in X
\end{equation}
and what can be said about the corresponding $g$'s. The answer is as follows:
\begin{lemma}\label{lemans} With $\Psi$, $\bar{x}_1$, $\bar{x}_2$ as above, $G\in\partial\Psi(\bar{x})$ satisfies {\rm (\ref{regularity})} if and only if the following two properties hold:
\begin{enumerate}
\item[(i)] $G$ is a ``certifying'' subgradient of $\Psi$ at $\bar{x}$, meaning that $\langle G,[0;x_2-\bar{x}_2]\rangle \geq0\,\,\forall x_2\in X_2[\bar{x}_1]$ (the latter relation indeed certifies that $\bar{x}_2$ is a minimizer of $\Psi(\bar{x}_1,x_2)$ over $x_2\in X_2[\bar{x}_1]$);
\item[(ii)] $g$ is a subgradient, taken at $\bar{x}_1$ along $X_1$, of the convex function
$$
\chi_G(x_1)=\min_{x_2\in X_2[x_1]} \langle G,[x_1;x_2]\rangle
$$
\end{enumerate}
\end{lemma}

It is easily seen that with $\Psi$, $\bar{x}=[\bar{x}_1;\bar{x}_2]$ as in Lemma \ref{lemans} (i.e., $\Psi$ is convex and Lipschitz continuous on $X$, $\bar{x}_1\in X_1$, and  $\bar{x}_2\in X_2[\bar{x}_1]$ minimizes $\Psi(\bar{x}_1,x_2)$ over $x_2\in X_2[\bar{x}_1]$) a certifying subgradient $G$ always exists; when $\Psi$ is differentiable at $\bar{x}$, one can take $G=\nabla_x\Psi(\bar{x})$. The function $\chi_G(\cdot)$, however, not necessary admits a subgradient at $\bar{x}_1$; when it does admit it, every $g\in\partial \chi_G(\bar{x}_1)$, satisfies (\ref{regularity}). In particular,
\begin{enumerate}
\item {[Direct Product case]} When $X= X_1\times X_2$, representing a certifying subgradient $G$ of $\Psi$, taken at  $[\bar{x}_1;\bar{x}_2\in\Argmin_{x_2\in X_2}\Psi(\bar{x}_1,x_2)]$, as $[g;h]$, we have
$$
\langle h,x_2-\bar{x}_2\rangle \geq0\,\,\forall x_2\in X_2,
$$
whence $\chi_G(x_1)=\langle g,x_1\rangle+\langle h,\bar{x}_2\rangle$, and thus $g$ is a subgradient of $\chi_G$ at $\bar{x}_1$. In particular, in the direct product case and when $\Psi$ is differentiable at $\bar{x}$, (\ref{regularity}) is met by $G=\nabla\Psi(\bar{x})$, $g=\nabla_{x_1}\Psi(\bar{x})$;
\item {[Polyhedral case]} When $X$ is a polyhedral set, for every certifying subgradient $G$ of $\Psi$ the function $\chi_G$ is polyhedrally representable with domain $X_1$ and as such has a subgradient at every point from $X_1$;
\item {[Interior case]} When $\bar{x}_1$ is a point from the relative interior of $X_1$, $\chi_G$ definitely  has a subgradient at $\bar{x}_1$.
\end{enumerate}

\subsection{Main Result, Saddle Point case}
\subsubsection{Preliminaries: execution protocols, accuracy certificates, residuals}
We start with outlining some simple concepts originating from \cite{NOR}. Let $W$ be a convex compact set in Euclidean space $\cW$, and $M(w):W\to \cW$ be a vector field on $W$. A {\sl $t$-step  execution protocol} associated with $M,W$ is a collection $\cI_t=\{w_i\in W, M(w_i):1\leq i\leq t\}$. A {\sl $t$-step accuracy certificate} is a $t$-dimensional probabilistic (i.e., with nonnegative entries summing up to 1) vector $\lambda$. Augmenting a $t$-step accuracy protocol $\cI_t$ by $t$-step accuracy certificate $\lambda$ gives rise to two entities:
\begin{equation}\label{twoent}
\begin{array}{rl}
\hbox{approximate solution:}&w^t=w^t(\cI_t,\lambda):=\sum_{i=1}^t \lambda_iw_i\in W;\\
\hbox{residual:}&\Res(\cI_t,\lambda_t|W)=\max\limits_{w\in W} \sum_{i=1}^t\lambda_i\langle M(w_i),w_i-w\rangle.\\
\end{array}
\end{equation}
When $W=U\times V$, where $U$ is a closed convex subset of Euclidean space $\cU$ and $V$ is a closed convex subset of Euclidean space $\cV$, and $M$ is vector field induced by  convex-concave function $\theta(u,v):U\times V\to\bR$, that is,
 $$M(u,v)=[M_u(u,v);M_v(u,v)]:U\times V\to \cU\times \cV\hbox{\ with\ } F_u(u,v)\in\partial_u\theta(u,v), \, F_v(u,v)\in\partial_v[-\theta(u,v)]$$ (such a field always is monotone),
 an execution protocol associated with $(M,W)$ will be called also {\sl  protocol associated with $\theta$, $U$, $V$}, or {\sl protocol associated with the saddle point problem} $$\min_{u\in U}\max_{v\in V}\theta(u,v).$$\par
The importance of these notions in our context stems from the following simple observation \cite{NOR}:
\begin{proposition}\label{lemNOR}
 Let $U$, $V$ be nonempty convex compact domains in Euclidean spaces $\cU$, $\cV$, $\theta(u,v):U\times V\to\bR$ be a convex-concave function,
 and $M$ be induced monotone vector field: $M(u,v)=[M_u(u,v);M_v(u,v)]:U\times V\to \cU\times \cV$ with $M_u(u,v)\in\partial_u\theta(u,v)$, $M_v(u,v)\in\partial_v[-\theta(u,v)]$. For a $t$-step execution protocol $\cI_t=\{w_i=[u_i;v_i]\in W:=U\times V,
 M_i=[M_u(u_i,v_i);M_v(u_i,v_i)]:1\leq i\leq t\}$ associated with $\theta,U,V$, and $t$-step accuracy certificate $\lambda$, it holds
 \begin{equation}\label{itholds}
 \epsilonsad(w^t(\cI_t,\lambda)|\theta,U,V)\leq\Res(\cI_t,\lambda|U\times V).
 \end{equation}
 \end{proposition}
\noindent Indeed, for $[u;v]\in U\times V$ we have
$$
\begin{array}{l}
\Res(\cI_t,\lambda|U\times V)\geq \sum_{i=1}^t\lambda_i\langle M_i,w_i-[u;v]\rangle=
\sum_{i=1}^t\lambda_i[\underbrace{\langle M_u(u_i,v_i),u_i-u\rangle}_{\geq \theta(u_i,v_i)-\theta(u,v_i)}
-\underbrace{\langle M_v(u_i,v_i),v_i-v\rangle}_{\leq \theta(u_i,v_i)-\theta(u_i,v)}]\\
\geq \sum_{i=1}^t\lambda_i[\theta(u_i,v)-\theta(u,v_i)]\geq \theta(u^t,v)-\theta(u,v^t),\\
\end{array}
$$
where the inequalities are due to the origin of $M$ and convexity-concavity of $\theta$. The resulting inequality holds true for all $[u;v]\in U\times V$, and (\ref{itholds}) follows.\qed
\subsubsection{Main Result}
\begin{proposition}\label{mlemma} In the situation and notation of sections \ref{s:sit} -- \ref{s:reg}, let $\phi$ be the primal convex-concave function induced by $\Phi$, and let
$$\cI_t=\{[x_{1,i};y_{1,i}]\in X_1\times Y_1,[\alpha_i:=\phi^\prime_{x_1}(x_{1,i},y_{1,i});\beta_i:=-\phi^\prime_{y_1}(x_{1,i},y_{1,i})]:1\leq i\leq t\}$$
be an execution protocol associated with $\phi$, $X_1,$ $Y_1$, where $\phi^\prime_{x_1}$, $\phi^\prime_{y_1}$ are regular sub- and supergradient fields associated with $\Phi$, $\phi$. Due to the origin of $\phi$, $\phi^\prime_{x_1}$, $\phi^\prime_{y_1}$, there exist $x_{2,i}\in X_2[x_{1,i}]$, $G_i\in \cX$, $y_{2,i}\in Y_2[y_{1,i}]$, and $H_i\in \cY$
 such that
\begin{equation}\label{eq3}
\begin{array}{lrcl}
(a)&G_i&\in&\partial_{x}\Phi(x_i:=[x_{1,i};x_{2,i}],y_i:=[y_{1,i};y_{2,i}]),\\
(b)&H_i&\in&\partial_{y}\left[-\Phi(x_i:=[x_{1,i};x_{2,i}],y_i:=[y_{1,i};y_{2,i}])\right],\\
(c)&\langle G_i,x-[x_{1,i};x_{2,i}]\rangle&\geq&\langle \phi_{x_1}^\prime(x_{1,i},y_{1,i}),x_1-x_{1,i}\rangle \forall x=[x_1;x_2]\in X,\\
(d)&\langle H_i,y-[y_{1,i};y_{2,i}]\rangle&\geq&\langle -\phi_{y_1}^\prime(x_{1,i},y_{1,i}),y_1-y_{1,i}\rangle \forall y=[y_1;y_2]\in Y,\\
\end{array}
\end{equation}
implying that
$$
\cJ_t=\left\{z_i=[x_i=[x_{1,i};x_{2,i}];y_i=[y_{1,i};y_{2,i}]],F_i=[G_i;H_i]:1\leq i\leq t\right\}
$$
is an execution protocol  associated with $\Phi$, $X$, $Y$. For every accuracy certificate $\lambda$ it holds
\begin{equation}\label{eq4}
\Res(\cJ_t,\lambda|X\times Y)\leq \Res(\cI_t,\lambda|X_1\times Y_1).
\end{equation}
As a result, given an accuracy certificate $\lambda$ and setting
$$
[x^t;y^t]=[[x_1^t;x_2^t];[y_1^t;y_2^t]]=\sum_{i=1}^t\lambda_i\left[[x_{1,i};x_{2,i}];[y_{1,i};y_{2,i}]\right],
$$
we ensure that
\begin{equation}\label{eq5}
\epsilonsad([x^t;y^t]|\Phi,X,Y)\leq \Res(\cI_t,\lambda|X_1\times Y_1),
\end{equation}
whence also, by Lemma \ref{lem1},
\begin{equation}\label{eq55}
\begin{array}{rcl}
\epsilonsad([x_1^t;y_1^t]|\phi,X_1,Y_1)&\leq& \Res(\cI_t,\lambda|X_1\times Y_1),\\
\epsilonsad([x_2^t;y_2^t]|\psi,X_2,Y_2)&\leq& \Res(\cI_t,\lambda|X_1\times Y_1),\\
\end{array}
\end{equation}
where $\psi$ is the dual function induced by $\Phi$.
\end{proposition}
{\bf Proof.} Let $z:=[[u_1;u_2];[v_1;v_2]]\in X\times Y$. Then
$$
\begin{array}{l}
\sum_{i=1}^t\lambda_i\langle F_i,z_i-z\rangle =\sum_{i=1}^t\lambda_i\bigg[\underbrace{\langle G_i,[x_{1,i};x_{2,i}]-[u_1;u_2]\rangle}_{\leq
\langle \phi_{x_1}^\prime(x_{1,i},y_{1,i}),x_{1,i}-u_1\rangle \hbox{\tiny\ by $(\ref{eq3}.c)$}}
 +
\underbrace{\langle H_i,[y_{1,i};y_{2,i}]-[v_1;v_2]\rangle}_{\leq \langle -\phi_{y_1}^\prime(x_{1,i},y_{1,i}),y_{1,i}-v_1\rangle \hbox{\tiny\ by $(\ref{eq3}.d)$}}\bigg] \\
\leq \sum_{i=1}^t\lambda_i\big[\langle \alpha_i,x_{1,i}-u_1\rangle +
\langle \beta_i,y_{1,i}-v_1\rangle\big]\leq \Res(\cI_t,\lambda|X_1\times Y_1), \\
\end{array}
$$
and (\ref{eq4}) follows. \qed
\subsection{Application: Solving bilinear  Saddle Point problems on domains represented by Linear Minimization Oracles}
\subsubsection{Situation}\label{sect:situ}
Let $W$ be a nonempty convex compact set in $\bR^N$, $Z$   be a nonempty convex compact set in $\bR^M$, and
let $\psi:W\times Z\to\bR$ be bilinear convex-concave function:
\begin{equation}\label{eq114}
\psi(w,z)=\langle w,p\rangle + \langle z,q\rangle + \langle z,Sw\rangle.
\end{equation}
Our goal is to solve the convex-concave SP problem
\begin{equation}\label{eq20}
 \min_{w\in W}\max_{z\in Z} \psi(w,z)
 \end{equation}
given by $\psi$, $W$, $Z$.
\subsubsection{Simple observation}
We intend to show that $\psi$ can be represented (in fact, in many ways) as the dual function induced by a bilinear convex-concave function $\Phi$; this is the key element of the outlined in section \ref{s:iccf} strategy for solving (\ref{eq20}).
\par
In the situation described in section \ref{sect:situ}, let $U\subset \bR^n$, $V\subset \bR^m$ be convex compact sets, and let $D\in\bR^{m\times N}$, $A\in\bR^{n\times M}$, $R\in\bR^{m\times n}$.  Consider bilinear (and thus convex-concave) function
\begin{equation}\label{eq114a}
\Phi(u,w;v,z)=\langle w,p+D^Tv\rangle + \langle z,q+A^Tu\rangle -\langle v,Ru\rangle:[U\times W]\times[V\times Z]\to\bR
\end{equation}
(the ``convex'' argument is $(u,w)$, the ``concave'' one is $(v,z)$). Assume that a pair of functions
\begin{equation}\label{eq107}
\begin{array}{l}
\bar{u}(w,z):W\times Z\to U,\\
\bar{v}(w,z):W\times Z\to V\\
\end{array}
\end{equation}
satisfies
\begin{equation}\label{eq115}
\begin{array}{lrcl}
\forall (w,z)\in W\times Z:&Dw&=&R\bar{u}(w,z)\\
\forall (w,z)\in W\times Z:&Az&=&R^T\bar{v}(w,z)\\
\end{array}
\end{equation}
Denoting $\bar{u}=\bar{u}(w,z)$, $\bar{v}=\bar{v}(w,z)$, we have
\begin{equation}\label{eq116}
\begin{array}{lrcl}
(a)&\langle w,D^T\bar{v}\rangle&=&\langle Dw,\bar{v}\rangle=\langle R\bar{u},\bar{v}\rangle\\
(b)&\langle z,A^T\bar{u}\rangle&=&\langle Az,\bar{u}\rangle=\langle \bar{u},R^T\bar{v}\rangle=\langle R\bar{u},\bar{v}\rangle.\\
\end{array}
\end{equation}
Thus,
$$
\begin{array}{rcl}
\nabla_u\Phi(\bar{u},w;\bar{v},z)&=&Az-R^T\bar{v}=0\\
\nabla_v\Phi(\bar{u},w;\bar{v},z)&=&Dw-R\bar{u}=0\\
\end{array}
$$
whence
$$
\begin{array}{rcl}
\bar{\psi}(w,z)&:=&\min_{u\in U}\max_{v\in V}\Phi(u,w;v,z)=\Phi(\bar{u}(w,z),w;\bar{v}(w,z),z)\\
&=&\langle w,p\rangle +\langle z,q\rangle +\langle Dw,\bar{v}(w,z)\rangle \hbox{\ [by (\ref{eq116})]}
\end{array}
$$
We have proved
\begin{lemma}\label{mainobs}
In the case of {\rm (\ref{eq107}), (\ref{eq115})}, assuming that
\begin{equation}\label{target}
\langle Dw,\bar{v}(w,z)\rangle =\langle z,Sw\rangle\,\,\forall (w\in W, z\in Z),
\end{equation}
$\psi$ is the dual convex-concave function induced by $\Phi$ and the domains $U\times W$, $V\times Z$.
\end{lemma}
Note that there are easy ways to ensure (\ref{eq115}) and (\ref{target}).
\paragraph{
Example 1.} Here $m=M$, $n=N$, and  $D=A^T=R=S$. Assuming $U\supset W$, $V\supset Z$ and setting  $\bar{u}(w,z)=w$, $\bar{v}(w,z)=z$, we ensure (\ref{eq107}), (\ref{eq115}) and (\ref{target}).
\paragraph{Example 2.} Let $S=A^TD$ with $A\in\bR^{K\times M}$, $D\in\bR^{K\times N}$. Setting $m=n=K$,
$R=I_K$, $\bar{u}(w,z)=Dw$, $\bar{v}(w,z)=Az$ and assuming that $U\supset DW$, $V\supset AZ$, we again ensure (\ref{eq107}), (\ref{eq115}) and  (\ref{target}).
\subsubsection{Implications}
Assume that (\ref{eq107}), (\ref{eq115}) and (\ref{target}) take place.
 Renaming the variables according to $x_1\equiv u$, $y_1\equiv v$, $x_2\equiv w$, $y_2\equiv z$ and setting $X_1=U$, $X_2=W$, $Y_1=V$, $Y_2=Z$, $X=X_1\times X_2=U\times W$, $Y=Y_1\times Y_2=V\times Z$, we find ourselves in the direct product case of the situation of section \ref{s:sit}, and
Lemma \ref{mainobs} says that the bilinear SP problem of interest (\ref{eq114}), (\ref{eq20})
  is the dual SP problem associated with the  bilinear master SP problem
 \begin{equation}\label{eq6}
 \min_{[u;w]\in U\times W}\max_{[v;z]\in V\times Z}\left[\Phi(u,w;v,z)=\langle w,p+D^Tv\rangle + \langle z,q+A^Tu\rangle -\langle Ru,v\rangle\right]
 \end{equation}
Since $\Phi$ is linear in $[w;z]$, the primal SP problem associated with (\ref{eq6}) is
$$
\min\limits_{u\equiv x_1\in U=X_1}\max\limits_{v\equiv y_1\in V=Y_1} \left[\phi(u,v)=\min\limits_{w\in W}\langle w,p+D^Tv\rangle + \max\limits_{z\in Z}\langle v,q+A^Tu\rangle -\langle Ru,v\rangle\right].
$$
Assuming that $W$, $Z$ allow for cheap Linear Minimization Oracles and defining $w_*(\cdot)$, $z_*(\cdot)$ according to
$$
w_*(\xi)\in\Argmin_{w\in W}\langle w,\xi\rangle,\,\,z_*(\eta)\in \Argmin_{z\in Z}\langle z,\eta\rangle,
$$
we have
\begin{equation}\label{eq22}
\begin{array}{rcl}
\phi(u,v)&=&\langle w_*(p+D^Tv),p+D^Tv\rangle +\langle z_*(-q-A^Tu),q+A^Tu\rangle - \langle Ru,v\rangle,\\
\phi^\prime_u(u,v)&:=&Az_*(-q-A^Tu)-R^Tv\in\partial_w\phi(u,v),\\
\phi^\prime_v(u,v)&:=&Dw_*(p+D^Tv)-Ru\in-\partial_v[-\phi(u,v)],\\
\end{array}
\end{equation}
that is, first order information on the primal SP problem
\begin{equation}\label{eq37}
\min\limits_{u\in U}\max\limits_{v\in V}\phi(u,v),
\end{equation}
is available. Note that since we are in the direct product case, $\phi^\prime_u$ and $\phi^\prime_v$ are regular sub- and supergradient fields associated with $\Phi$, $\phi$.
\par
Now let $\cI_t=\{[u_i;v_i]\in U\times V,[\gamma_i:=\phi^\prime_u(u_i,v_i);\delta_i:=-\phi^\prime_v(u_i,v_i)]:1\leq i\leq t\}$ be an execution protocol generated by a First Order algorithm as applied to the primal SP problem (\ref{eq37}), and let
$$
\begin{array}{l}
w_i=w_*(p+D^Tv_i),z_i=z_*(-q-A^Tu_i),\\
\alpha_i=\nabla_w\Phi(u_i,w_i;v_i,z_i)=p+D^Tv_i,\\
\beta_i=-\nabla_z\Phi(u_i,w_i;v_i,z_i)=-q-A^Tu_i,\\
\end{array}
$$
so that
$$
\cJ_t=\bigg\{[[u_i;w_i];[v_i;z_i]],[\underbrace{[\gamma_i;\alpha_i]}_{\hbox{\tiny\rm$\nabla_{[u;w]}\Phi(u_i,w_i;v_i,z_i)$}};
\underbrace{[\delta_i;\beta_i]}_{
\hbox{\tiny\rm$-\nabla_{[v;z]}\Phi(u_i,w_i;v_i,z_i)$}}]:1\leq i\leq t\bigg\}
$$
is an execution protocol associated with the SP problem (\ref{eq6}). By Proposition \ref{mlemma}, for any accuracy certificate $\lambda$ it holds
\begin{equation}\label{eq30}
\Res(\cJ_t,\lambda|U\times W\times V\times Z)\leq \Res(\cI_t,\lambda|U\times V)
\end{equation}
whence, setting
\begin{equation}\label{eq31}
[[u^t;w^t];[v^t;z^t]]=\sum_{i=1}^t\lambda_i[[u_i;w_i];[v_i;z_i]]
\end{equation}
and invoking Proposition \ref{lemNOR} with $\Phi$ in the role of $\theta$,
\begin{equation}\label{eq32}
 \epsilonsad([[u^t;w^t];[v^t;z^t]]|\Phi,\underbrace{X_1\times X_2}_{U\times W},\underbrace{Y_1\times Y_2}_{V\times Z})\leq \Res(\cI_t,\lambda|U\times V)
\end{equation}
whence, by Lemma \ref{lem1},
\begin{equation}\label{eq33}
\epsilonsad([w^t;z^t]|\psi,W,Z)\leq \Res(\cI_t,\lambda|U\times V).
\end{equation}
We have arrived at the following
\begin{proposition}\label{theSP} In the situation of section \ref{sect:situ}, let
{\rm (\ref{eq107}), (\ref{eq115})} and {\rm (\ref{target})} take place. Then applying to the primal SP problem {\rm (\ref{eq37})} First Order algorithm $\cB$ with accuracy certificates, $t$-step execution protocol $\cI_t$ and accuracy certificate $\lambda^t$ generated by $\cB$ yield straightforwardly a feasible solution to the SP problem of interest {\rm (\ref{eq114}) -- (\ref{eq20})} of the  $\epsilonsad$-inaccuracy $\leq \Res(\cI_t,\lambda^t|U\times V)$.
 \end{proposition}
 Note also that when the constructions from Examples 1, 2 are used, there is a significant freedom in selecting the domain $U\times V$ of the primal problem (we only require $U$, $V$ to be convex compact sets ``large enough'' to ensure the inclusions mentioned in Examples), so that there is no difficulty to enforce $U$, $V$ to be proximal friendly. As a result, we can take as $\cB$ a proximal First Order method, for example, Non-Euclidean Restricted Memory Level algorithm with certificates (cf. \cite{DualSubgradient}) or Mirror Descent (cf. \cite{FTVI}). The efficiency estimates of these algorithms as given in \cite{DualSubgradient,FTVI} imply  that the resulting procedure for solving the SP of interest (\ref{eq114}) -- (\ref{eq20}) admits non-asymptotic  $O(1/\sqrt{t})$ rate of convergence, with explicitly computable factors hidden in $O(\cdot)$. The resulting complexity bound is completely similar to the one achievable with the machinery of Fenchel-type representations  \cite{DualSubgradient,FTVI}.
\par We are about to consider a special case where the $O(1/\sqrt{t})$ complexity admits a significant improvement.
\subsection{Matrix Game case}
Let $S\in\bR^{M\times N}$ admit representation
 $$S=A^TD$$
 with $A\in\bR^{K\times M}$ and $D\in\bR^{K\times N}$. Let also $W=\Delta_N=\{w\in\bR^N_+: \sum_iw_i=1\}$,  $Z=\Delta_M$. Our goal is to solve matrix game
\begin{equation}\label{mgame}
\min_{w\in W}\max_{z\in Z}\left[\psi(w,z)=\langle z,Sw\rangle=\langle Az,Dw\rangle \right].
\end{equation}
Let $U$, $V$ be convex compact sets such that
\begin{equation}\label{suchthat}
V\supset AZ,\,\,U\supset DW,
\end{equation}
and let us set
$$
\begin{array}{rcl}
\Phi(u,w;v,z)&=&\langle u,Az\rangle+\langle v,Dw\rangle  -\langle u,v\rangle\\
\bar{u}:=\bar{u}(w,z)&=&Dw\\
\bar{v}:=\bar{v}(w,z)&=&Az\\
\end{array}
$$                                                                                                                                     implying that
$$
\begin{array}{rcl}
\nabla_u\Phi(\bar{u},w;\bar{v},z)&=&Az-\bar{v}=0,\\
\nabla_v\Phi(\bar{u},w;\bar{v},z)&=&Dw-\bar{u}=0,\\
\Phi(\bar{u},w;\bar{v},z)&=&\langle \bar{u},Az\rangle+\langle \bar{v},Dw\rangle -\langle \bar{u},\bar{v}\rangle=
\langle Dw,Az\rangle+\langle Az,Dw\rangle-\langle Dw,Az\rangle  \\
&=&\langle z,A^TDw\rangle=\psi(w,z).\\
\end{array}
$$
\def\Min{{\hbox{\rm Min}}}
\def\Max{{\hbox{\rm Max}}}
It is immediately seen that the function $\psi$ from (\ref{mgame}) is nothing but the dual convex-concave function associated with $\Phi$ (cf. Example 2), while the primal function is
\begin{equation}\label{eqphi}
\phi(u,v)=\Max(A^Tu)+\Min(D^Tv)-\langle u,v\rangle;
\end{equation}
here $\Min(p)$ and $\Max(p)$ stand for the smallest and the largest entries in vector $p$.
Applying the strategy outlined in section \ref{s:iccf}, we can solve the problem of interest (\ref{mgame}) applying to the primal
SP problem
\begin{equation}\label{pgame}
\min_{u\in U}\max_{v\in V} \left[\phi(u,v)=\Min(D^Tv)+\Max(A^Tu)-\langle u,v\rangle\right]
\end{equation}
an algorithm with accuracy certificates and using the machinery outlined in previous sections to convert the resulting execution protocols and certificates
into approximate solutions to the problem of interest (\ref{mgame}).
\par
We intend to consider a special case when the outlined approach allows to reduce a huge, but simple, matrix game (\ref{mgame}) to a small SP problem (\ref{pgame}) --  so small that it can be solved to high accuracy by a cutting plane method (e.g., the Ellipsoid algorithm). This is the case when the matrices $A$, $D$ in (\ref{mgame}) are {\sl simple}.
\subsubsection{Simple  matrices}
Given a $K\times L$ matrix $B$, we call $B$ {\sl simple} if, given $x\in\bR^K$, it is easy to identify the columns $\overline{B}[x]$, $\underline{B}[x]$ of $B$ making the maximal, resp. the minimal, inner product with $x$.
\par
 When matrices $A$, $D$ in (\ref{mgame}) are simple, the first order information for the cost function $\phi$ in the primal SP problem (\ref{pgame}) is easy to get. Besides, all we need from the convex compact sets $U$, $V$ participating in (\ref{pgame}) is to be large enough to ensure that $U\supset DW$ and $V\supset AZ$, which allows to make $U$ and $V$ simple, e.g., Euclidean balls. Finally, when the design dimension $2K$ of (\ref{pgame}) is small, we have at our disposal a multitude of linearly converging, with the converging ratio depending solely on $K$, methods for solving (\ref{pgame}), including the Ellipsoid algorithm with certificates presented in \cite{NOR}.  We are about to demonstrate that the outlined situation indeed takes place in some meaningful applications.
\subsubsection{Example: Knapsack generated matrices}\label{sect:knapsack}
\footnote{The construction to follow can be easily extended from ``knapsack generated'' matrices to more general ``Dynamic Programming generated'' ones, see section \ref{DP} in Appendix.} Assume that we are given {\sl knapsack data}, namely,
\begin{itemize}
\item positive integer {\sl horizon} $m$,
\item nonnegative integer {\sl bounds} $\bar{p}_s$, $1\leq s\leq m$,
\item positive integer {\sl costs} $h_s$, $1\leq s\leq m$, and positive integer {\sl budget} $H$,  and
\item {\sl output functions} $f_s(\cdot):\{0,1,...,\bar{p}_s\}\to\bR^{r_s}$, $1\leq s\leq m$.
\end{itemize}
Given the outlined data, consider the set $\cP$ of all integer vectors $p=[p_1;...;p_m]\in\bR^m$ satisfying the following restrictions:
$$
\begin{array}{ll}
0\leq p_s\leq\overline{p}_s,\,1\leq s\leq m&\hbox{[range restriction]}\\
\sum_{s=1}^mh_sp_s\leq H&\hbox{[budget restriction]}\\
\end{array}
$$
and the matrix $B$ of the size $K\times \Card(\cP)$, $K=\sum_{s=1}^m r_s$, defined as follows: the columns of $B$ are indexed by vectors $p=[p_1;...;p_s]\in\cP$, and the column indexed by $p$ is the vector
$$B_p=[f_1(p_1);...;f_m(p_m)].$$
Note that assuming $m$, $\overline{p}_s$, $r_s$ moderate, matrix $B$ is simple -- given $x\in\bR^K$, it is easy to find $\overline{B}[x]$ and $\underline{B}[x]$ by Dynamic Programming.
\begin{quote}
Indeed, to identify $\overline{B}[x]$, $x=[x_1;...;x_m]\in\bR^{r_1}\times...\times\bR^{r_m}$ (identification of $\underline{B}[x]$ is completely similar), it suffices to run for $s=m,m-1,...1$ the backward Bellman recurrence
$$
\left.
\begin{array}{rcl}
U_s(h)&=&\max\limits_{r\in\bZ}\left\{U_{s+1}(h-h_s r)+\langle f_s(r),x_s\rangle: 0\leq r\leq \overline{p}_s, 0\leq h-h_sr\right\}\\
A_s(h)&\in&\Argmax\limits_{r\in\bZ}\left\{U_{s+1}(h-h_s r)+\langle f_s(r),x_s\rangle: 0\leq r\leq \overline{p}_s, 0\leq h-h_sr\right\}\\
\end{array}\right\}, h=0,1,...,H,
$$
with $U_{m+1}(\cdot)\equiv0$,  and then to recover one by one the entries $p_s$ in the index $p\in\cP$ of $\overline{B}[x]$ from the forward Bellman recurrence
$$
\begin{array}{l}
H_1=H,p_1=A_1(H_1);\\
H_{s+1}=H_s-h_sp_s,p_{s+1}=A_{s+1}(H_{s+1}),1\leq s<m.
\end{array}
$$
\end{quote}
\par
\subsubsection{Illustration: Attacker vs. Defender.}\label{Blotto} The ``covering story'' we intend to consider is as follows\footnote{This story is a variation of what is called ``Colonel Blotto Game'' in Game Theory, see, e.g., \cite{Bellman,Robertson} and references therein.}. Attacker and Defender are preparing for a  conflict to take place on $m$ battlefields. A pure strategy of Attacker is a vector $a=[a_1;...;a_m]$, where nonnegative integer $a_s$, $1\leq s\leq m$,  is the number of attacking units to be created and deployed at battlefield $s$; the only restrictions on $a$, aside of nonnegativity and integrality, are the bounds $a_s\leq\overline{a}_s$, $1\leq s\leq m$,  and the budget constraint $\sum_{s=1}^mh_{sA} a_s\leq H_A$ with positive integer $h_{sA}$ and $H_A$.  Similarly, a pure strategy of Defender is  a vector $d=[d_1;...;d_m]$, where nonnegative integer $d_s$ is the number of defending units to be created and deployed at battlefield $s$, and the only restrictions on $d$, aside of nonnegativity and integrality, are the bounds $d_s\leq\overline{d}_s$, $1\leq s\leq m$,  and the budget constraint $\sum_{s=1}^mh_{sD} d_s\leq H_D$ with positive integer $h_{sD}$ and $H_D$. The total loss of Defender (the total gain of Attacker), the pure strategies of the players being $a$ and $d$, is
$$
S_{a,d}=\sum_{s=1}^m [\Omega^s]_{a_s,d_s},
$$
with given $(\overline{a}_s+1)\times(\overline{d}_s+1)$ matrices $\Omega^s$. Our goal is to solve in mixed strategies
the matrix game where Defender seeks to minimize his total loss, and Attacker seeks to maximize it.
\par
Denoting by $\cA$ and $\cD$ the sets of pure strategies of Attacker, resp., Defender, representing
$$
\Omega^s=\sum_{i=1}^{r_s}f^{is}[g^{is}]^T,\,f^{is}=[f^{is}_0;...;f^{is}_{\overline{a}_s}],\,\,
g^{is}=[g^{is}_0;...;g^{is}_{\overline{d}_s}],\,\,r_s=\rank(\Omega^s),
$$
and setting
$$
\begin{array}{rcl}
K&=&\sum_{s=1}^m r_s,\\
A_a&=&[[f^{1,1}_{a_1};f^{2,1}_{a_1};...;f^{r_1,1}_{a_1}];[f^{1,2}_{a_2};f^{2,2}_{a_2};...;f^{r_2,2}_{a_2}];...;
[f^{1,m}_{a_m};f^{2,m}_{a_m};...;f^{r_m,m}_{a_m}]]\in\bR^K,a\in\cA,\\
D_d&=&[[g^{1,1}_{d_1};g^{2,1}_{d_1};...;g^{r_1,1}_{d_1}];[g^{1,2}_{d_2};g^{2,2}_{d_2};...;g^{r_2,2}_{d_2}];...;
[g^{1,m}_{d_m};g^{2,m}_{d_m};...;g^{r_m,m}_{d_m}]]\in\bR^K,d\in\cD,\\
\end{array}
$$
we end up with $K\times M$, $M=\Card(\cA)$, knapsack-generated matrix $A$ with columns $A_a$, $a\in \cA$, and $K\times N$, $N=\Card(\cD)$, knapsack-generated matrix $D$ with columns $D_d$, $d\in\cD$, such that
$$
S:=[S_{a,d}]_{{a\in\cA\atop d\in\cD}}=A^TD.
$$
As a result,  solving the Attacker vs. Defender game in  mixed strategies reduces to solving SP  problem (\ref{mgame}) with knapsack-generated (and  thus simple) matrices $A$, $D$ and thus can be reduced to convex-concave SP (\ref{pgame}) on the product of two $K$-dimensional convex compact sets. Note that in the situation in question the design dimension $2K$ of (\ref{pgame}) will, typically, be rather small (few tens or at most few hundreds), while the design dimensions $M$, $N$ of the matrix game of interest (\ref{mgame}) can be huge.
 \paragraph{Numerical illustration.} With the data (quite reasonable in terms of the ``Attacker vs. Defender'' game)
$$
m=8,\, h_{sA}=h_{sD}=1,1\leq s\leq m,\, H_A=H_D=64=\overline{d}_s=\overline{a}_s,\,1\leq s\leq m
$$
and rank 1 matrices $\Omega_s$, $1\leq s\leq m$, the design dimensions of the problem of interest (\ref{mgame}) are as huge as
$$
\dim w=\dim z=97,082,021,465
$$
while the sizes of problem (\ref{pgame}) are just
$$
\dim u=\dim v=8,
$$
 and thus (\ref{pgame}) can be easily solved to high accuracy by the Ellipsoid method. In the numerical experiment we are about to report\footnote{for implementation details, see section \ref{ellM}}, the outlined approach allowed to solve (\ref{mgame}) within $\epsilonsad$-inaccuracy as small as 5.0e-9 in just 1537 steps of the Ellipsoid algorithm (110.0 sec on a medium quality laptop).  This performance is quite promising, especially
when taking into account huge -- nearly 10$^{11}$ -- sizes of the matrix game of interest (\ref{mgame}).

\section{From  Saddle Point problems to Variational Inequalities with Monotone Operators}\label{sect:monotoneVI}
In what follows, we extend the decomposition approach (developed so far for convex-concave SP problems) to Variational Inequalities (VI's) with monotone operators, with the primary goal to handle VI's with affine monotone operators on LMO-represented domains.
\subsection{Decomposition of Variational Inequalities with monotone operators}\label{sect3.1}
\subsubsection{Situation}\label{sect:situation}
Let $\cX$, $\cH$ be Euclidean spaces, $\Theta\subset \cX\times \cH$ be convex compact set, $\Xi$ be the projection of $\Theta$ onto $\cX$, and $H$ be the projection of  $\Theta$ onto $\cH$. Given $\xi\in \Xi$, $\eta\in H$, we set
$$
H_\xi=\{\eta:[\xi;\eta]\in \Theta\},\,\,\Xi_\eta=\{\xi\in \Xi: [\xi;\eta]\in \Theta\}.
$$
We denote a point from $\cX\times \cH$ as $\theta=[\xi;\eta]$ with $\xi\in \cX$, $\eta\in \cH$. Let, further, $$\Phi(\xi,\eta)=[\Phi_\xi(\xi,\eta);\Phi_\eta (\xi,\eta)]:\Theta\to \cX \times \cH $$ be a continuous monotone vector field.

\subsubsection{Induced vector fields}
Let $\xi \in \Xi$, and let $\overline{\eta}=\overline{\eta}(\xi)$ be a somehow selected, as a function of $\xi \in \Xi$, strong solution to the VI given by $(H_{\xi},\Phi_\eta (\xi,\eta))$, that is,
\begin{equation}\label{etazero}
\overline{\eta}(\xi)\in H_{\xi}\ \&\ \langle \Phi_\eta (\xi,\overline{\eta}(\xi)),\eta -\overline{\eta}(\xi)\rangle \geq0\,\,\forall \eta \in H_{\xi}.
\end{equation}
Let us call $\Phi$ (more precisely, the pair $(\Phi,\overline{\eta}(\cdot))$) {\sl $\eta $-regular}, if for every $\xi \in \Xi$, there exists $\Psi=\Psi(\xi)
\in \cX $ such that
\begin{equation}\label{vregularity}
\langle \Psi(\xi),\xi'-\xi \rangle \leq \langle \Phi(\xi,\overline{\eta}(\xi)),[\xi';\eta']-[\xi;\overline{\eta}(\xi)]\rangle \,\forall [\xi';\eta']\in \Theta.
\end{equation}
Similarly, let $\overline{\xi}(\eta)$ be a somehow selected, as a function of $\eta\in H$, strong solution to the VI given by $(\Xi_\eta,\Phi_\eta(\xi,\eta))$, that is,
\begin{equation}\label{xizero}
\overline{\xi}(\eta)\in \Xi_{\eta}\ \&\ \langle \Phi_\xi (\overline{\xi}(\eta),\eta),\xi -\overline{\xi}(\eta)\rangle \geq0\,\,\forall \xi \in \Xi_{\eta}.
\end{equation}
Let us call $(\Phi,\overline{\xi}(\cdot))$ {\sl $\xi $-regular}, if for every $\eta \in H$ there exists $\Gamma=\Gamma(\eta)\in \cH $ such that
\begin{equation}\label{uregularity}
\langle \Gamma(\eta),\eta'-\eta \rangle \leq \langle \Phi(\overline{\xi}(\eta),\eta),[\xi';\eta']-[\overline{\xi}(\eta);\eta]\rangle \,\forall [\xi';\eta']\in \Theta.
\end{equation}
When $(\Phi,\overline{\eta})$ is $\eta$-regular, we refer to the above $\Psi(\cdot)$ as to a {\sl primal} vector field induced by $\Phi$  \footnote{``a primal'' instead of ``the primal'' reflects the fact that $\Psi$ is not uniquely defined by $\Phi$ -- it is defined by $\Phi$ {\sl and} $\overline{\eta}$ {\sl and} by how the values of $\Psi$ are selected when (\ref{vregularity}) does not specify these values uniquely.}, and when $(\Phi,\overline{\xi})$ is $\xi$-regular, we refer to the above $\Gamma(\cdot)$ as to a {\sl dual} vector field induced by $\Phi$.
\paragraph{Example: Direct product case.} This is the case where $\Theta=\Xi\times H$. In this situation, setting $\Psi(\xi)=\Phi_\xi (\xi,\overline{\eta}(\xi))$, we have for $[\xi';\eta']\in \Theta$:
$$
\langle \Phi(\xi,\overline{\eta}(\xi)),[\xi';\eta']-[\xi;\overline{\eta}(\xi)]\rangle=\underbrace{\langle \Phi_\xi (\xi,\overline{\eta}(\xi)),\xi'-\xi \rangle}_{=\langle \Psi(\xi),\xi'-\xi\rangle} +\underbrace{\langle \Phi_\eta (\xi,\overline{\eta}(\xi)),\eta'-\overline{\eta}(\xi)\rangle}_{\geq 0\,\,\forall \eta'\in H_\xi =H}\geq
\langle \Psi(\xi),\xi'-\xi \rangle,
$$
that is, $(\Phi,\overline{\eta}(\cdot))$ is $\eta $-regular, with $\Psi(\xi)=\Phi_\xi (\xi,\overline{\eta}(\xi))$. Setting $\Gamma(\eta)=
\Phi_\eta(\overline{\xi}(\eta),\eta)$, we get by similar argument
$$
\langle \Phi(\overline{\xi}(\eta)),[\xi';\eta']-[\overline{\xi};\eta]\rangle\geq \langle \Gamma(\eta),\eta'-\eta\rangle, \eta,\eta'\in H,
$$
that is, $(\Phi,\overline{\xi}(\cdot))$ is $\xi$-regular, with $\Gamma(\eta)=
\Phi_\eta(\overline{\xi}(\eta),\eta).$
\subsubsection{Main Result, Variational Inequality case}
\subsubsection{Preliminaries}
Recall that the (Minty's) variational inequality $\VI(M,W)$ associated with a convex compact subset $W$ of Euclidean space $\cW$ and a vector field $M:W\to \cW$ is
$$
\hbox{find\ }w\in W: \langle M(w'),w'-w\rangle \geq0\,\,\forall w'\in W\eqno{\VI(M,W)}
$$
$w$ satisfying the latter condition is called a {\sl weak} solution to the VI. A natural measure of inaccuracy for an approximate solution $w\in W$ to $\VI(M,W)$ is the {\sl dual gap function}
$$
\epsilonVI(w|M,W)=\sup_{w'\in W} \langle M(w'),w-w'\rangle;
$$
weak solutions to the VI are exactly the points of $W$ where this (clearly nonnegative everywhere on $W$) function is zero.
\par
In the sequel we utilize the following simple fact originating from \cite{NOR}:
\begin{proposition}\label{lemNORA} Let $M$ be monotone on $W$, let $\cI_t=\{w_i\in W,M(w_i):1\leq i\leq t\}$ be a $t$-step execution protocol associated with $(M,W)$, $\lambda$ be a $t$-step accuracy certificate, and $w^t=\sum_{i=1}^t\lambda_iw_i$ be the associated approximate solution. Then
$$
\epsilonVI(w^t|M,W)\leq\Res(\cI_t,\lambda|W).
$$
\end{proposition}
\noindent
Indeed, we have
$$
\begin{array}{l}
\Res(\cI_t,\lambda|W)=\sup_{w'\in W}\left[\sum_{i=1}^t\lambda_i\langle M(w_i),w_i-w'\rangle\right]\\
\geq\sup_{w'\in W}\left[ \sum_{i=1}^t\lambda_i\langle M(w'),w_i-w'\rangle\right] \hbox{\ [since $M$ is monotone]}\\
=\sup_{w'\in W}\langle M(w'), w^t-w'\rangle =\epsilonVI(w^t|M,W).\\
\end{array}\eqno{\hbox{\qed}}
$$
\subsubsection{Main result}
\begin{proposition}\label{lemMonOp} In the situation of section \ref{sect:situation}, let $(\Phi,\overline{\eta}(\cdot))$ be $\eta $-regular. Then
\par
{\rm (i)} Primal  vector field
$\Psi(\xi)$ induced by $(\Phi,\overline{\eta}(\cdot))$ is monotone on $\Xi$. Moreover, whenever $\cI_t=\{\xi_i\in \Xi,\Psi(\xi_i):1\leq i\leq t\}$ and $\cJ_t=\{\theta_i:=[\xi_i;\overline{\eta}(\xi_i)],\Phi(\theta_i):1\leq i\leq t\}$ and $\lambda$ is a $t$-step  accuracy certificate, it holds
\begin{equation}\label{itholdsVI}
\epsilonVI(\sum_{i=1}^t\lambda_i\theta_i|\Phi,\Theta)\leq \Res(\cJ_t,\lambda|\Theta)\leq\Res(\cI_t,\lambda|\Xi).
\end{equation}
\par
{\rm (ii)} Let $(\Phi,\overline{\xi})$ be $\xi$-regular, and let $\Gamma$ be the induced dual vector field. Whenever $\widehat{\theta}=[\widehat{\xi};\widehat{\eta}]\in \Theta$, we have
\begin{equation}\label{follows}
\epsilonVI(\widehat{\eta}|\Gamma,H)\leq\epsilonVI(\widehat{\theta}|\Phi,\Theta).
\end{equation}
\end{proposition}

\subsection{Implications}\label{sect3.2}
In the situation of section \ref{sect:situation}, assume that for properly selected $\overline{\eta}(\cdot)$, $\overline{\xi}(\cdot)$, $(\Phi,\overline{\eta}(\cdot))$ is $\eta$-regular, and $(\Phi,\overline{\xi}(\cdot))$ is $\xi$-regular, induced primal and dual vector fields being $\Psi$ and $\Gamma$. In order to solve the {\sl dual VI} $\VI(\Gamma,H)$, we can apply to the {\sl primal VI} $\VI(\Psi,\Xi)$ an algorithm with accuracy certificates; by Proposition \ref{lemMonOp}.i, resulting $t$-step execution protocol $\cI_t=\{\xi_i,\Psi(\xi_i):1\leq i\leq t\}$ and accuracy certificate $\lambda$ generate an execution protocol $\cJ_t=\{\theta_i:=[\xi_i;\overline{\eta}(\xi_i)],\Phi(\theta_i):1\leq i\leq t\}$ such that
$$
\Res(\cJ_t,\lambda|\Theta)\leq \Res(\cI_t,\lambda|\Xi),
$$
whence, by Proposition \ref{lemNORA}, for the approximate solution
$$
\theta^t=[\xi^t,\eta^t]:=\sum_{i=1}^T\lambda_i\theta_i=\sum_{i=1}^t\lambda_i[\xi_i;\overline{\eta}(\xi_i)]
$$
it holds
$$\epsilonVI(\theta^t|\Phi,\Theta)\leq \Res(\cI_t,\lambda|\Xi).
$$
Invoking Proposition \ref{lemMonOp}.ii, we conclude that $\eta^t$ is a feasible solution to the dual VI $\VI(\Gamma,H)$, and
\begin{equation}\label{ofinterest}
\epsilonVI(\eta^t|\Gamma,H)\leq \Res(\cI_t,\lambda|\Xi).
\end{equation}
We are about to present two examples well suited for the just outlined approach.
\subsubsection{Solving affine monotone VI on LMO-represented domain}
Let $H$ be a convex compact set in $\cH=\bR^N$, and let $H$ be equipped with an LMO. Assume that we want to solve the VI $\VI(F,H)$, where
$$
F(\eta)=S\eta+s
$$
is an affine monotone operator (so that $S+S^T\succeq0$). Let us set
$\cX=\cH$, select $\Xi$ as a proximal-friendly convex compact set containing $H$, and set $\Theta=\Xi\times H$,
$$
\Phi(\xi,\eta)=\underbrace{\left[\begin{array}{c|c}S^T &-S^T\\
\hline
S&\cr\end{array}\right]}_{\cS}\left[\begin{array}{c}\xi\cr\eta\cr\end{array}\right]+\left[\begin{array}{c}0\cr
s\cr\end{array}\right].
$$
We have
$$
\cS+\cS^T=\left[\begin{array}{c|c}S+S^T&\cr\hline
&\cr \end{array}\right]\succeq0,
$$
so that $\Phi$ is an affine monotone operator with
$$
\begin{array}{rcl}
\Phi_\xi(\xi,\eta)&=&S^T\xi-S^T\eta\\
\Phi_\eta(\xi,\eta)&=&S\xi+s\\
\end{array}
$$
Setting $\overline{\xi}(\eta)=\eta$, we ensure that $\overline{\xi}(\eta)\in \Xi$ when $\eta\in H$ and $\Phi_\xi(\overline{\xi}(\eta),\eta)=0$,   implying (\ref{xizero}). Since we are in the direct product case, we can set
$\Gamma(\eta)=\Phi_\eta(\overline{\xi}(\eta),\eta)=S\eta+s=F(\eta)$; thus, $\VI(\Gamma,H)$ is our initial VI of interest.
On the other hand, setting
$$
\overline{\eta}(\xi)\in\Argmin_{\eta\in H} \langle S\xi+s,\eta\rangle,
$$
we ensure (\ref{etazero}). Since we are in the direct product case, we can set
$$
\Psi(\xi)=\Phi_\xi(\xi,\overline{\eta}(\xi))=S^T[\xi-\overline{\eta}(\xi)];
$$
note that the values of $\Psi$ can be straightforwardly computed via calls to the LMO representing $H$.
We can now solve $\VI(\Psi,\Xi)$ by a proximal algorithm $\cB$ with accuracy certificates and recover, as explained above, approximate solution to the VI of interest $\VI(F,H)$. With the Non-Euclidean Restricted Memory Level method with certificates \cite{DualSubgradient} or Mirror Descent with certificates (see, e.g., \cite{FTVI}), the approach results in non-asymptotical $O(1/\sqrt{t})$-converging algorithm for solving the VI of interest, with explicitly computable factors hidden in $O(\cdot)$. This complexity bound, completely similar to the one obtained in \cite{FTVI}, seems to be the best known under the circumstances.
\subsubsection{Solving skew-symmetric VI on LMO-represented domain}\label{sect:skewsymm}
Let $H$ be an LMO-represented convex compact domain in $\cH=\bR^N$, and assume that we want to solve  $\VI(F,H)$, where
$$
F(\eta)=2Q^TP\eta+f:\cH\to \cH
$$
with $K\times N$ matrices $P,Q$ such that the matrix $Q^TP$ is skew-symmetric:
\begin{equation}\label{eq888}
Q^TP+P^TQ=0.
\end{equation}
Let $\cX=\bR^K\times\bR^K$, let $\Xi_1$, $\Xi_2$ be two convex compact sets in $\bR^K$ such that
\begin{equation}\label{suchthatNew}
Q H\subset \Xi_1, -PH\subset \Xi_2.
\end{equation}
Let us set $\Xi=\Xi_1\times\Xi_2$, and let
$$
\Phi(\xi=[\xi_1;\xi_2],\eta)=
\left[\begin{array}{c|c|c}&I_K&P\cr\hline
-I_K&&Q\cr\hline
-P^T&-Q^T&\cr\end{array}\right]\left[\begin{array}{c}\xi_1\cr \xi_2\cr
\eta\cr\end{array}\right]+\left[\begin{array}{c}0\cr 0\cr f\cr\end{array}\right].
$$
Note that $\Phi$ is monotone and affine. Setting
$$
\overline{\xi}(\eta)=[Q\eta;-P\eta]
$$
and invoking (\ref{suchthatNew}), we ensure (\ref{xizero}); since we are in the direct product case, we can take, as the dual induced vector field,
$$
\Gamma(\eta)=\Phi_{\eta}(\overline{\xi}(\eta),\eta)=-P^T(Q\eta)-Q^T(-P\eta)+f=[Q^TP-P^TQ]\eta+f\underbrace{=}_{\hbox{\tiny by (\ref{eq888})}} 2Q^TP\eta+f=F(\eta),
$$
so that the dual VI $\VI(\Gamma,H)$ is our VI of interest.
\par
On the other hand, setting
$$
\overline{\eta}(\xi=[\xi_1;\xi_2])\in\Argmin_{\eta\in H} \langle f-P^T\xi_1-Q^T\xi_2,\eta\rangle,
$$
we ensure (\ref{etazero}). Since we are in the direct product case, we can define primal vector field as
$$
\Psi(\xi=[\xi_1;\xi_2])=\Phi_\xi([\xi_1;\xi_2],\overline{\eta}([\xi_1;\xi_2]))=\left[\begin{array}{c}\xi_2+P\overline{\eta}(\xi)\cr
-\xi_1+Q\overline{\eta}(\xi)\cr\end{array}\right].
$$
Note that LMO for $H$ allows to compute the values of $\Psi$, and that $\Xi$ can be selected to be proximal-friendly.
We can now solve $\VI(\Psi,\Xi)$ by a proximal algorithm $\cB$ with accuracy certificates and recover, as explained above, approximate solution to the VI of interest $\VI(F,H)$. When the design dimension $\dim \Xi$ of the primal VI is small, other choices of $\cB$, like the Ellipsoid algorithm, are possible, and in this case we can end up with linearly converging, with the converging ratio depending solely on $\dim \Xi$, algorithm for solving the VI of interest. We are about to give a related example, which can be considered as multi-player version of the ``Attacker vs. Defender'' game.
\paragraph{Example: Nash Equilibrium with pairwise interactions.} Consider the situation as follows: there are
 \begin{itemize}
 \item $L\geq2$ players, $\ell$-th of them selecting a mixed strategy $w_\ell$ from probabilistic simplex $\Delta_{N_\ell}$ of dimension $N_\ell$,
 \item {\sl encoding matrices} $D_\ell$ of sizes $m_\ell\times N_\ell$, and {\sl loss matrices} $M^{\ell\ell'}$ of sizes $m_\ell\times m_{\ell'}$ such that
$$
M^{\ell\ell}=0,M^{\ell\ell'}=-[M^{\ell'\ell}]^T,\,1\leq\ell,\ell'\leq L.
$$
\item
The loss of $\ell$-th player depends on mixed strategies of the players according to
$$
\cL_\ell(\eta:=[w_1;...;w_L])=\sum_{\ell'=1}^L w_\ell^TE^{\ell\ell'}w_{\ell'},\,\,E^{\ell\ell'}=D_\ell^TM^{\ell\ell'}D_{\ell'} +\langle g_\ell,\eta\rangle.
$$
In other words, every pair of distinct players $\ell,\ell'$ are playing matrix game with matrix $M^{\ell\ell'}$, and the loss of player $\ell$, up to a linear in $[w_1;...;w_L]$ function, is the sum, over the pairwise games he is playing, of his losses in these games, the ``coupling constraints'' being expressed by the requirement that every player uses the same mixed strategy in all pairwise games he is playing.
\end{itemize}
We have described {\sl convex} Nash Equilibrium problem, meaning that for every $\ell$, $\cL_\ell(w_1,...,w_L)$ is convex (in fact, linear) in $w_\ell$, is jointly concave (in fact, linear) in $w^\ell:=(w_1,...,w_{\ell-1},w_{\ell+1},...,w_L)$, and $\sum_{\ell=1}^L\cL_\ell(\eta)$ is the linear function $\langle g,\eta\rangle$, $g=\sum_\ell g_\ell$, and thus is convex. It is known (see, e.g., \cite{NOR}) that Nash  Equilibria in convex Nash problem are exactly the weak solutions to the VI given by monotone operator
$$
F(\eta:=[w_1;...;w_L])=[\nabla_{w_1}\cL_1(\eta);...;\nabla_{w_L}\cL_L(\eta)]
$$
on the domain
$$
H=\Delta_{N_1}\times...\times \Delta_{N_L}.
$$
Let us set
$$
Q={1\over 2}\left[\begin{array}{cccc}D_1&&&\\
&D_2&&\\
&&\ddots&\\
&&&D_L\\
\end{array}\right],\,\,
P=\left[\begin{array}{cccc} M^{1,1}D_1&M^{1,2}D_2&...&M^{1,L}D_L\\
M^{2,1}D_1&M^{2,2}D_2&...&M^{2,L}D_L\\
\vdots&\vdots&\ddots&...\\
M^{L,1}D_1&M^{L,2}D_2&...&M^{L,L}D_L\\
\end{array}\right].\\
$$
Then
$$
Q^TP={1\over 2}\left[\begin{array}{cccc}D_1^TM^{1,1}D_1&D_1^TM^{1,2}D_2&...&D_1^TM^{1,L}D_L\\
D_2^TM^{2,1}D_1&D_2^TM^{2,2}D_2&...&D_2^TM^{1,L}D_L\\
\vdots&\vdots&\ddots&...\\
D_L^TM^{L,1}D_1&D_L^TM^{L,2}D_2&...&D_L^TM^{L,L}D_L\\
\end{array}\right]
$$
so that $Q^TP$ is skew-symmetric due to $M^{\ell\ell'}=-[M^{\ell'\ell}]^T$. Besides this, we clearly have
$$
F(\eta:=[w_1;...;w_L])=2Q^TP\eta+f, \,\,f=[\nabla_{w_1}\langle g_1,[w_1;...;w_L]\rangle;...;\nabla_{w_L}\langle g_L,[w_1;...;w_L]\rangle].
$$
Observe that if $D_1,...,D_L$ are simple, so are $Q$ and $P$.
 \begin{quote}
 Indeed, for $Q$ this is evident: to find the column of $Q$ which makes the largest inner product with $x=[x_1;...;x_L]$, $\dim x_\ell=m_\ell$, it suffices to find, for every $\ell$, the column of $D_\ell$ which makes the maximal inner product with $x_\ell$, and then to select the maximal of the resulting $L$ inner products and the corresponding to this maximum column of $Q$. To maximize the inner product of the same $x$ with columns of $P$, note that
$$
x^TP=\big[\underbrace{[{\sum}_{\ell=1}^Lx_\ell^TM^{\ell,1}]}_{y_1^T}D_1,...,\underbrace{[{\sum}_{\ell=1}^Lx_\ell^TM^{\ell,L}]}_{y_L^T}D_L\big],
$$ so that to maximize the inner product of $x$ and the columns of $P$ means to find, for every $\ell$,  the column of $D_\ell$ which makes the maximal inner product with $y_\ell$, and then to select the maximal of the resulting $L$ inner products and the corresponding to this maximum column of $P$.
\end{quote}
We see that if $D_\ell$ are simple, we can use the approach from  section \ref{sect:skewsymm} to approximate the solution to the VI generated by $F$ on $H$. Note that in the case in question the dual gap function $\epsilonVI(\eta|F,H)$ admits a transparent interpretation in terms of the Nash Equilibrium problem we are solving: for $\eta=[w_1;...;w_L]\in H$, we have
\begin{equation}\label{wehaveew}
\epsilonVI(\eta|F,H)\geq\epsilonNash(\eta):=\sum_{\ell=1}^L \left[\cL_\ell(\eta)-\min_{w_\ell^\prime\in\Delta_{N_\ell}}\cL_\ell(w_1,...,w_{\ell-1},w_\ell^\prime,w_{\ell+1},...,w_L)\right],
\end{equation}
and the right hand side here is the sum, over the players, of the (nonnegative) incentives for a player $\ell$ to deviate from his strategy $w_\ell$ to another mixed strategy when all other players stick to their strategies as given by $\eta$. Thus, small $\epsilonVI([w_1;...;w_L]|\cdot,\cdot)$ means small incentives for the players to deviate from mixed strategies $w_\ell$.
\begin{quote}
Verification of (\ref{wehaveew}) is immediate: denoting $f_\ell=\nabla_{w_\ell}\langle g_\ell,w\rangle$, by definition of $\epsilonVI$ we have for every $\eta'=[w_1^\prime;...;w_L^\prime]\in H$:
$$
\begin{array}{l}
\epsilonVI(\eta|F,H)\geq \langle F(\eta'),\eta-\eta'\rangle
 =\sum_{\ell} \langle \nabla_{w_\ell}\cL_\ell(\eta'),w_\ell-w_\ell^\prime\rangle\\ =\sum_{\ell}\langle f_\ell,w_\ell-w_\ell^\prime\rangle +\sum_{\ell,\ell'}\langle D_\ell^TM^{\ell\ell'}D_{\ell'}w_{\ell'}^\prime,w_\ell-w_\ell^\prime\rangle \\
=\sum_{\ell}\langle f_\ell,w_\ell-w_\ell^\prime\rangle +\sum_{\ell,\ell'}\langle D_\ell^TM^{\ell\ell'}D_{\ell'}w_{\ell'},w_\ell-w_\ell^\prime\rangle \\
\multicolumn{1}{r}{\hbox{[since $\sum_{\ell,\ell'}\langle D_\ell^TM^{\ell\ell'}D_{\ell'}z_{\ell'},z_\ell\rangle=0$ due to $M^{\ell\ell'}=-[M^{\ell'\ell}]^T$]}}\\
=\sum_{\ell}\langle \nabla_{w_\ell}\cL(\eta),w_\ell-w_\ell^\prime\rangle = \sum_{\ell}[\cL_\ell(\eta)-\cL_\ell(w_1,...,w_{\ell-1},w_\ell^\prime,w_{\ell+1},...,w_L)]\\
\multicolumn{1}{r}{\hbox{[since $\cL_\ell$ is affine in $w_\ell$]}}\\
\end{array}
$$
and (\ref{wehaveew}) follows.
\end{quote}
\subsection{Relation to \cite{FTVI}}\label{sect:dominance}
Here we demonstrate that the decomposition approach to solving VI's with monotone operators on LMO-represented domains cover the approach, based on Fenchel-type representations, developed in \cite{FTVI}. Specifically, let $H$ be a compact convex set in Euclidean space $\cH=\bR^N$, $G(\cdot)$ be a monotone vector field on $H$, and $\eta\mapsto Ax+a$ be an affine mapping from $\cH$ to Euclidean space $\cX=\bR^M$. Given a convex compact set $\Xi\subset \cX$, let us set
\begin{equation}\label{setup}
\Theta=\Xi\times H, \,\,\Phi(\xi,\eta)=\left[\Phi_\xi(\xi,\eta):=A\eta + a;\Phi_\eta(\xi,\eta):=G(\eta)-A^T\xi\right]:
\Theta\to \cX\times \cH,
\end{equation}
so that $\Phi$ clearly is a monotone vector field on $\Theta$. Assume that  $\overline{\eta}(\xi):\Xi\to H$ is a somehow selected strong solution to $\VI(\Phi_\eta(\xi,\cdot),H)$:
\begin{equation}\label{repr}
\forall \xi\in\Xi: \overline{\eta}(\xi)\in H\ \&\ \underbrace{\langle G(\overline{\eta}(\xi))-A^T\xi,\eta-\overline{\eta}(\xi)\rangle}_{=\langle
\Phi_\eta(\xi,\overline{\eta}(\xi)),\eta-\overline{\eta}(\xi)\rangle}
\geq0\,\forall \eta\in H;
\end{equation}
(cf. (\ref{etazero}));
note that required $\overline{\eta}(\xi)$ definitely exists, provided that $G(\cdot)$ is continuous and monotone. Let us also define $\overline{\xi}(\eta)$ as a selection of the point-to-set mapping $\eta\mapsto \Argmin\limits_{\xi\in \Xi}\langle A\eta+a,\xi\rangle$, so that
\begin{equation}\label{repra}
\forall \eta\in H: \overline{\xi}(\eta)\in \Xi\ \&\
\underbrace{\langle A\eta+a,\xi-\overline{\xi}(\eta)\rangle}_{=\langle \Phi_\xi(\overline{\xi}(\eta),\eta),\xi-\overline{\xi}(\eta)\rangle}\geq0,\forall \xi\in\Xi
\end{equation}
(cf. (\ref{xizero})).
\par
Observe that with the just defined $\Xi$, $H$, $\Theta$, $\Phi$, $\overline{\eta}(\cdot)$, $\overline{\xi}(\cdot)$ we are in the direct product case of the situation described in section \ref{sect:situation}.
Since we are in the direct product case, $(\Phi,\overline{\eta}(\cdot))$ is $\eta$-regular, and we can take, as the induced primal vector field associated with $(\Phi,\overline{\eta}(\cdot))$, the vector field
\begin{equation}\label{Psi}
\Psi(\xi)=A\overline{\eta}(\xi)+a=\Phi_\xi(\xi,\overline{\eta}(\xi)):\Xi\to \cX,
\end{equation}
and as the induced dual vector field, the  field
\[
\Gamma(\eta)=G(\eta)-A^T\overline{\xi}(\eta)=\Phi_\eta(\overline{\xi}(\eta),\eta):H\to \cX.
\]
\par
Note that in terms of \cite{FTVI}, relations (\ref{Psi}) and (\ref{repr}),  modulo notation, form what in the reference is called a {\sl Fenchel-type representation} (F.-t.r.) {\sl of vector field $\Psi:\Xi\to \cX$, the data of the representation being $\cH$, $A$, $a$, $\overline{\eta}(\cdot)$, $G(\cdot)$, $H$}. On a closer inspection, every F.-t.r. of a given monotone vector field $\Psi:\Xi\to \cX$
can be obtained in this fashion  from some setup of the form (\ref{setup}). \par
Assume now that $\Xi$ is LMO-representable, and we have at our disposal $G$-oracle which, given on input $\eta\in H$, returns $G(\eta)$.
 This oracle combines with LMO for $\Xi$ to induce a procedure which, given on input $\eta\in H$, returns $\Gamma(\eta)$. As a result, we can apply the decomposition machinery presented in sections \ref{sect3.1},
\ref{sect3.2} to reduce solving $\VI(\Psi,\Xi)$ to processing $\VI(\Gamma,H)$ by an algorithm with accuracy certificates. It can be easily seen by inspection that this reduction recovers constructions and results presented in \cite[sections 1 -- 4]{FTVI}. The bottom line is that the  developed in sections \ref{sect3.1}, \ref{sect3.2} decomposition-based approach to solving VI's with monotone operators on LMO-represented domains essentially covers the  developed in \cite{FTVI} approach based on Fenchel-type representations of monotone vector fields\footnote{``covers'' instead of ``is equivalent'' stems from the fact that the scope of decomposition is not restricted to the setups of the form of (\ref{setup})).}.
\appendix
\section{Appendix}
\subsection{Proof of Lemma \ref{lem1}}
It suffices to prove the $\phi$-related statements. Lipschitz continuity of $\phi$ in the direct product case is evident. Further, the function $\theta(x_1,x_2;y_1)=\max\limits_{y_2\in Y_2[y_1]}\Phi(x_1,x_2;y_1,y_2)$ is convex and Lipschitz continuous in $x=[x_1;x_2]\in X$ for every $y_1\in Y_1$, whence
$
\phi(x_1,y_1)=\min\limits_{x_2\in X_2[x_1]}\theta(x_1,x_2;y_1)$ is convex and lower semicontinuous in $x_1\in X_1$ (note that $X$ is compact). On the other hand, $\phi(x_1,y_1)=\max\limits_{y_2\in Y_2[y_1]}\min\limits_{x_2\in X_2[x_1]}\Phi(x_1,x_2;y_1,y_2)=\max\limits_{y_2\in Y_2[y_1]} \left[\chi(x_1;y_1,y_2):=\min\limits_{x_2\in X_2[x_1]}\Phi(x_1,x_2;y_1,y_2)\right]$, so that $\chi(x_1;y_1,y_2)$ is concave and Lipschitz continuous in $y=[y_1;y_2]\in Y$
for every $x_1\in X_1$, whence $$\phi(x_1,y_1)=\max\limits_{y_2\in Y_2[y_1]}\chi(x_1;y_1,y_2)$$
is concave and upper semicontinuous in $y_1\in Y_1$ (note that $Y$ is compact).
\par
Next, we have
$$
\begin{array}{l}
\SV(\phi,X_1,X_2)=\inf\limits_{x_1\in X_1}\left[\sup\limits_{y_1\in Y_1}\left[\sup\limits_{y_2:[y_1;y_2]\in Y}\inf\limits_{x_2:
[x_1;x_2]\in X}\Phi(x_1,x_2;y_1,y_2)\right]\right]\\
=
\inf\limits_{x_1\in X_1}\left[\sup\limits_{[y_1;y_2]\in Y}\inf\limits_{x_2:[x_1;x_2]\in X}\Phi(x_1,x_2;y_1,y_2)\right]\\
=\inf\limits_{x_1\in X_1}\left[\inf\limits_{x_2:[x_1;x_2]\in X}\sup\limits_{[y_1;y_2]\in Y}\Phi(x_1,x_2;y_1,y_2)\right]\hbox{\ [by Sion-Kakutani Theorem]}\\
=\inf\limits_{[x_1;x_2]\in X}\sup\limits_{[y_1;y_2]\in Y}\Phi(x_1,x_2;y_1,y_2)=\SV(\Phi,X,Y),\\
\end{array}
$$
as required in (\ref{eq9}). Finally, let $\bar{x}=[\bar{x}_1;\bar{x}_2]\in X$ and $\bar{y}=[\bar{y}_1;\bar{y}_2]\in Y$. We have
$$
\begin{array}{rcl}
\overline{\phi}(\bar{x}_1)-\SV(\phi,X_1,Y_1)&=&\overline{\phi}(\bar{x}_1)-\SV(\Phi,X,Y)\hbox{\ [by (\ref{eq9})]}\\
&=&\sup\limits_{y_1\in Y_1}\phi(\bar{x}_1,y_1)-\SV(\Phi,X,Y)\\
&=&\sup\limits_{y_1\in Y_1}\sup\limits_{y_2:[y_1;y_2]\in Y}\inf\limits_{x_2:[\bar{x}_1;x_2]\in X}\Phi(\bar{x}_1,x_2;y_1,y_2)-\SV(\Phi,X,Y)\\
&=&\sup\limits_{[y_1;y_2]\in Y}\inf\limits_{x_2:[\bar{x}_1;x_2]\in X}\Phi(\bar{x}_1,x_2;y_1,y_2)-\SV(\Phi,X,Y)\\
&=&\inf\limits_{x_2:[\bar{x}_1;x_2]\in X}\sup\limits_{y=[y_1;y_2]\in Y}\Phi(\bar{x}_1,x_2;y)-\SV(\Phi,X,Y)\\
&\leq&\sup\limits_{y=[y_1;y_2]\in Y}\Phi(\bar{x}_1,\bar{x}_2;y)-\SV(\Phi,X,Y)\\
&=&\overline{\Phi}(\bar{x})-\SV(\Phi,X,Y)\\
\end{array}
$$
and
$$
\begin{array}{rcl}
\SV(\phi,X_1,Y_1)-\underline{\phi}(\bar{y}_1)&=&\SV(\Phi,X,Y)-\underline{\phi}(\bar{y}_1)\hbox{\ [by (\ref{eq9})]}\\
&=&\SV(\Phi,X,Y)-\inf\limits_{x_1\in X_1}\phi(x_1,\bar{y}_1)\\
&=&\SV(\Phi,X,Y)-\inf\limits_{x_1\in X_1}\left[\inf\limits_{x_2:[x_1;x_2]\in X}\sup\limits_{y_2:[\bar{y}_1;y_2]\in Y}\Phi(x_1,x_2;\bar{y}_1,y_2)\right]\\
&=&\SV(\Phi,X,Y)-\inf\limits_{x=[x_1;x_2]\in X}\sup\limits_{y_2:[\bar{y}_1;y_2]\in Y}\Phi(x;\bar{y}_1,y_2)\\
&\leq& \SV(\Phi,X,Y)-\inf\limits_{x=[x_1;x_2]\in X}\Phi(x;\bar{y}_1,\bar{y}_2)\\
&=&\SV(\Phi,X,Y)-\underline{\Phi}(\bar{y}).\\
\end{array}
$$
We conclude that
$$
\begin{array}{l}
\epsilonsad([\bar{x}_1;\bar{y}_1]|\phi,X_1,Y_1)=\left[\overline{\phi}(\bar{x}_1)-\SV(\phi,X_1,Y_1)\right]+\left[\SV(\phi,X_1,Y_1)-\underline{\phi}(\bar{y}_1)\right]\\
\leq \left[\overline{\Phi}(\bar{x})-\SV(\Phi,X,Y)\right]+\left[\SV(\Phi,X,Y)-\underline{\Phi}(\bar{y})\right]=
\epsilonsad([\bar{x};\bar{y}]|\Phi,X,Y),\\
\end{array}
$$
as claimed in (\ref{eq10}). \qed
\subsection{Proof of Lemma \ref{lem2}}
For $x_1\in X_1$ we have
$$
\begin{array}{l}
\phi(x_1;\bar{y}_1) =\min\limits_{x_2:[x_1;x_2]\in X}\max\limits_{y_2:[\bar{y}_1;y_2]\in Y}\Phi(x_1,x_2;\bar{y}_1,y_2)\geq \min\limits_{x_2:[x_1;x_2]\in X}\Phi(x_1,x_2;\bar{y}_1,\bar{y}_2)\\
\geq \min\limits_{x_2:[x_1;x_2]\in X}\big[\underbrace{\Phi(\bar{x};\bar{y})}_{\phi(\bar{x}_1;\bar{y}_1))}+\langle G,[x_1;x_2]-[\bar{x}_1;\bar{x}_2]\big]\rangle\hbox{\ [since $\Phi(x;\bar{y})$ is convex and $G\in\partial_x\Phi(\bar{x};\bar{y})$]}\\
\geq  \phi(\bar{x}_1;\bar{y}_1)+\langle g,x_1-\bar{x}_1\rangle\hbox{\ [by definition of $g,G$]},\\
\end{array}
$$
as claimed in $(a)$. ``Symmetric'' reasoning justifies $(b)$. \qed
\subsection{Proof of Lemma \ref{lemans}}
Assume that (\ref{regularity}) holds true. Then $G$ clearly is certifying, implying that
$$
\chi_G(\bar{x}_1)=\langle G,[\bar{x}_1;\bar{x}_2]\rangle,
$$
and therefore (\ref{regularity}) reads
$$
\langle G,[x_1;x_2]\rangle \geq \chi_G(\bar{x}_1)+\langle g,x_1-\bar{x}_1\rangle \,\,\forall x=[x_1;x_2]\in X,
$$
where, taking minimum in the left hand side over $x_2\in X_2[x_1]$,
$$
\chi_G(x_1)\geq \chi_G(\bar{x}_1)+\langle g,x_1-\bar{x}_1\rangle\,\,\forall x_1\in X_1,
$$
as claimed in (ii).\par
Now assume that (i) and (ii) hold true. By (i), $\chi_G(\bar{x}_1)=\langle G,[\bar{x}_1;\bar{x}_2]\rangle$, and by (ii) combined with the definition of $\chi_G$,
$$
\forall x=[x_1;x_2]\in X: \langle G,[x_1;x_2]\rangle \geq \chi_G(x_1)\geq \chi_G(\bar{x}_1)+\langle g,x_1-\bar{x}_1\rangle = \langle G,\bar{x}\rangle +
\langle g,x_1-\bar{x}_1\rangle,$$ implying (\ref{regularity}). \qed

\subsection{Dynamic Programming generated simple matrices}\label{DP}
Consider the situation as follows. There exists an evolving in time system $\cS$, with state $\xi_s$ at time $s=1,2,...,m$ belonging to a given finite nonempty set $\Xi_s$. Further, every pair $(\xi,s)$ with $s\in \{1,...,m\}$, $\xi\in \Xi_s$ is associated with nonempty finite {\sl set of actions} $A^s_\xi$, and we set
$$
\cS_s=\{(\xi,a):\xi\in \Xi_s,a\in A^s_\xi\}.
$$
 Further, for every $s$, $1\leq s< m$, a {\sl transition mapping} $\pi_{s}(\xi,a):\cS_s\to \Xi_{s+1}$ is given. Finally, we are given vector-valued functions (''outputs'') $\chi_s:\cS_s\to\bR^{r_s}$.\par
  A {\sl trajectory} of $\cS$ is a sequence $\{(\xi_s,a_s):1\leq s\leq m\}$ such that $(\xi_s,a_s)\in \cS_s$ for $1\leq s\leq m$ and
$$
 \xi_{s+1}=\pi_{s}(\xi_s,a_s),\,1\leq s<m.
 $$
The {\sl output} of a trajectory $\tau=\{(\xi_s,a_s):1\leq s\leq m\}$ is the block-vector $$\chi[\tau]=[\chi_1(\xi_1,a_1);...;\chi_m(\xi_m,a_m)].$$ We can associate with $\cS$ the matrix $D=D[\cS]$ with $K=r_1+...+r_m$ rows and with columns indexed by the trajectories of $\cS$; specifically, the column indexed by a trajectory $\tau$ is $\chi[\tau]$.
\begin{quote}
For example, knapsack generated matrix $D$ associated with knapsack data from section \ref{sect:knapsack} is of the form $D[\cS]$ with system $\cS$ as follows:
\begin{itemize}
\item $\Xi_s$, $s=1,...,m$, is the set of nonnegative integers which are $\leq H$;
\item $A^s_\xi$ is the set of nonnegative integers $a$ such that $a\leq\bar{p}_s$ and $\xi-h_sp_s\geq0$;
\item the transition mappings are $\pi_{s}(\xi,a)=\xi-ah_s$;
\item the outputs are $\chi_s(\xi,a)=f_s(a)$, $1\leq s\leq m$.
\end{itemize}
In the notation from section \ref{sect:knapsack}, vectors $[p_1;...;p_m]\in\cP$ are exactly the sequences of actions $a_1,...,a_m$ stemming from the trajectories of the just defined system $\cS$. 
\end{quote}
Observe that matrix $D=D[\cS]$ is simple, provided the cardinalities of $\Xi_s$ and $A^s_\xi$ are reasonable.
Indeed, given $x=[x_1;...;x_m]\in\bR^{n}=\bR^{r_1}\times...\times\bR^{r_m}$, we can identify $\overline{D}[x]$ by Dynamic Programming, running first the backward Bellman recurrence
$$
\left.\begin{array}{rcl}
U_{s}(\xi)&=&\max\limits_{a\in A^s_\xi}\left\{x_s^T\chi_s(\xi,a)+U_{s+1}(\pi_s(\xi,a))\right\}\\
A_s(\xi)&=&\Argmax\limits_{a\in A^s_\xi}\left\{x_s^T\chi_s(\xi,a)+U_{s+1}(\pi_s(\xi,a))\right\}\\
\end{array}\right\}, \xi\in\Xi_s,\, s=m,m-1,...,1
$$
(where $U_{m+1}(\cdot)\equiv 0$),
and then identifying the (trajectory indexing the) column of $D$ corresponding to $\overline{D}[x]$ by running the forward Bellman recurrence
$$
\begin{array}{rcl}
\xi_1&\in&\Argmax_{\xi\in \Xi_1} U_1(\xi)\Rightarrow a_1\in A_1(\xi_1)\Rightarrow...\\
\Rightarrow \xi_{s+1}&=&\pi_s(\xi_s,a_s)\Rightarrow a_{s+1}\in A_{s+1}(\xi_{s+1})\Rightarrow ...\\
\end{array}, s=1,2,...,m-1.
$$
\subsection{Attacker vs. Defender via Ellipsoid algorithm}\label{ellM}

In our implementation,
\begin{enumerate}
\item Relation (\ref{suchthatNew}) is ensured by specifying $U$, $V$ as centered at the origin Euclidean balls of radius $R$,  where $R$ is an upper bound on the Euclidean norms of the columns in $D$ and in $A$ (such a bound can be easily obtained from the knapsack data specifying the matrices $D$, $A$).
\item We processed the monotone vector field associated with the primal SP problem (\ref{pgame}), that is, the field
$$
F(u,v)=[F_u(u,v)=\overline{A}[u]-v;F_v(u,v)=u-\underline{D}[v]]
$$
by Ellipsoid algorithm with accuracy certificates from \cite{NOR}. For $\tau=1,2,...,$ the algorithm generates  {\sl search points}
$[u_\tau;v_\tau]\in\bR^K\times\bR^K$, with $[u_1;v_1]=0$, along with execution protocols $\cI^\tau=\{[u_i;v_i],F(u_i,v_i):i\in I_\tau\}$, where $I_\tau=\{i\leq\tau:[u_i;v_i]\in U\times V\}$, augmented by accuracy certificates $\lambda^\tau=\{\lambda^\tau_i\geq0:i\in I_\tau\}$ such that $\sum_{i\in I_\tau}\lambda^\tau_i=1$. From the results of \cite{NOR} it follows that for every $\epsilon>0$ it holds
\begin{equation}\label{ellMet}
\tau\geq N(\epsilon):= O(1)K^2\ln\left(2{R+\epsilon\over\epsilon}\right)\Rightarrow \Res(\cI^\tau,\lambda^\tau|U\times V)\leq\epsilon.
\end{equation}
\item When computing  $F(u_i,v_i)$ (this computation takes place only at {\sl productive} steps -- those with $[u_i;v_i]\in U\times V$), we get, as a byproduct, the columns $A^i=\overline{A}[u_i]$ and $D^i=\underline{D}[v_i]$ of matrices $A$, $D$, along with the indexes $a^i$, $d^i$  of these columns (recall that these indexes are pure strategies of Attacker and Defender and thus, according to the construction of $A$, $D$, are collections of $m$ nonnegative integers).  In our implementation, we stored these columns, same as their indexes and the corresponding search points $[u_i;v_i]$. As is immediately seen, in the case in question the approximate solution $[w^\tau;z^\tau]$ to the SP problem of interest (\ref{mgame}) induced by execution protocol $\cI^\tau$ and accuracy certificate $\lambda^\tau$ is comprised of two sparse vectors
\begin{equation}\label{eqwz}
w^\tau=\sum_{i\in I_\tau}\lambda^\tau_i\delta^D_{d^i},\,\,z^\tau=\sum_{i\in I_\tau}\lambda^\tau_i\delta^A_{a^i}
\end{equation}
where $\delta^D_d$ is the ``$d$-th basic orth'' in the simplex $\Delta_N$ of probabilistic vectors with entries indexed by pure strategies of Defender, and similarly for $\delta^A_a$. Thus, we have no difficulties with representing our approximate solutions\footnote{Note that applying Caratheodory theorem, we could further ``compress'' the representations of approximate solutions -- make these solutions convex combinations of at most $K+1$ of $\delta^D_{d^i}$'s and  $\delta^A_{a^i}$'s.}, in spite of their huge ambient dimension.
\end{enumerate}
According to our general theory and (\ref{ellMet}), the number of steps needed to get an $\epsilon$-solution $[w;z]$ to the problem of interest (i.e., a feasible solution with $\epsilonsad([w;z]|\psi,W,Z)\leq\epsilon)$ does not exceed $N(\epsilon)$, with computational effort per step dominated by the necessity to identify $\overline{A}[u_i]$, $\underline{D}[v_i]$ by Dynamic Programming.
\par
In fact, we used the outlined scheme with two straightforward modifications.
\begin{itemize}
\item First, instead of building the accuracy certificates $\lambda^\tau$ according to the rules from \cite{NOR}, we used the best, given execution protocols $\cI^\tau$, accuracy certificates
by solving the convex program
\begin{equation}\label{cvxproblem}
\min_\lambda \left\{\Res(\cI^\tau,\lambda|U\times V):=\max_{y\in U\times V}\sum_{i\in I_\tau} \lambda_i
\langle F(u_i,v_i),[u_i;v_i]-y\rangle:\lambda_i\geq0,\sum_{i\in I_\tau}\lambda_i=1\right\}
\end{equation}
In our implementation, this problem was solved from time to time, specifically, once per $4K^2$ steps. Note that with $U$, $V$ being Euclidean balls, (\ref{cvxproblem})  is well within the scope of Matlab Convex Programming solver {\tt CVX} \cite{cvx}.
\item Second, given current approximate solution (\ref{eqwz}) to the problem of interest, we can compute its saddle point inaccuracy exactly instead of upper-bounding it by $\Res(\cI^\tau,\lambda^\tau|U\times V)$. Indeed, it is immediately seen that
$$
\epsilonsad([w^\tau;z^\tau]|\psi,W,Z)=\Max(A^T[\sum_{i\in I_\tau}\lambda^\tau_iD^i])-\Min(D^T[\sum_{i\in I_\tau}\lambda^\tau_iA^i]).
$$
In our implementation, we performed this computation each time when a new accuracy certificate was computed, and terminated the solution process when the saddle point inaccuracy became less than a given threshold (1.e-4).
\end{itemize}
\subsection{Proof of Proposition \ref{lemMonOp}}
(i): Let $\xi_1,\xi_2\in \Xi$, and let $\eta_1=\overline{\eta}(\xi_1)$, $\eta_2=\overline{\eta}(\xi_2)$. By (\ref{vregularity}) we have
$$
\begin{array}{rcl}
\langle \Psi(\xi_2),\xi_2-\xi_1)&\geq&\langle \Phi(\xi_2,\eta_2),[\xi_2-\xi_1;\eta_2-\eta_1]\rangle\\
\langle \Psi(\xi_1),\xi_1-\xi_2)&\geq&\langle \Phi(\xi_1,\eta_1),[\xi_1-\xi_2;\eta_1-\eta_2]\rangle\\
\end{array}
$$
Summing inequalities up, we get
$$
\langle \Psi(\xi_2)-\Psi(\xi_1),\xi_2-\xi_1\rangle\geq \langle \Phi(\xi_2,\eta_2)-\Phi(\xi_1,\eta_1),[\xi_2-\xi_1;\eta_2-\eta_1]\rangle \geq 0,
$$
so that $\Psi$ is monotone.
\par
Further, the first inequality in (\ref{itholdsVI}) is due to Proposition \ref{lemNORA}. To prove the second inequality in (\ref{itholdsVI}), let $\cI_t=\{\xi_i\in \Xi,\Psi(\xi_i):1\leq i\leq t\}$, $\cJ_t=\{\theta_i:=[\xi_i;\overline{\eta}(\xi_i)],\Phi(\theta_i):1\leq i\leq t\}$, and let $\lambda$ be $t$-step accuracy certificate. We have
$$
\begin{array}{l}
\theta=[\xi;\eta ]\in \Theta \Rightarrow\\
\sum_{i=1}^t\lambda_i\langle \Phi(\theta_i),\theta_i-\theta\rangle\leq \sum_{i=1}^t\lambda_i \langle \Psi(\xi_i),\xi_i-\xi \rangle
 \hbox{\ [see (\ref{vregularity})]}\\
 \leq \Res(\cI_t,\lambda|\Xi)\\
\Rightarrow
\Res(\cJ_t,\lambda|\Theta)=\sup_{\theta=[\xi;\eta ]\in \Theta} \sum_{i=1}^t\lambda_i\langle \Phi(\theta_i),\theta_i-\theta\rangle\leq \Res(\cI_t,\lambda|\Xi).
\end{array}
$$
(i) is proved.
\par
(ii): Let $\eta \in H$. Invoking (\ref{uregularity}), we have
$$
\langle \Gamma(\eta),\widehat{\eta}-\eta \rangle \leq \langle \Phi(\overline{\xi}(\eta),\eta), [\widehat{\xi};\widehat{\eta}]-[\overline{\xi}(\eta);\eta]\rangle \leq \epsilonVI(\widehat{\theta}|\Phi,\Theta),
$$
and (\ref{follows})  follows. \qed

\end{document}